\documentclass{article}

\usepackage{microtype}
\usepackage{graphicx}
\usepackage{subfigure}
\usepackage{booktabs} 

\usepackage{amssymb, amsmath}
\usepackage{xcolor}
\usepackage{algpseudocode}

\usepackage{hyperref}

\algblockdefx{MRepeat}{EndRepeat}{\textbf{repeat}}{}
\algnotext{EndRepeat}

\usepackage[accepted]{icml2020}

\icmltitlerunning{Responsive Safety in RL by PID Lagrangian Methods}

\begin{document}

\twocolumn[
\icmltitle{Responsive Safety in Reinforcement Learning by PID Lagrangian Methods}



\icmlsetsymbol{equal}{*}

\begin{icmlauthorlist}
\icmlauthor{Adam Stooke}{ucb,openai}
\icmlauthor{Joshua Achiam}{ucb,openai}
\icmlauthor{Pieter Abbeel}{ucb}
\end{icmlauthorlist}

\icmlaffiliation{ucb}{University of California, Berkeley}
\icmlaffiliation{openai}{OpenAI}

\icmlcorrespondingauthor{Adam Stooke}{adam.stooke@berkeley.edu}

\icmlkeywords{Reinforcement Learning, Safety, Safe Reinforcement Learning, Constrained Optimization, Lagrangian Method}

\vskip 0.3in
]



\printAffiliationsAndNotice{}  

\begin{abstract}

Lagrangian methods are widely used algorithms for constrained optimization problems, but their learning dynamics exhibit oscillations and overshoot which, when applied to safe reinforcement learning, leads to constraint-violating behavior during agent training.  We address this shortcoming by proposing a novel Lagrange multiplier update method that utilizes derivatives of the constraint function.  We take a controls perspective, wherein the traditional Lagrange multiplier update behaves as \emph{integral} control; our terms introduce \emph{proportional} and \emph{derivative} control, achieving favorable learning dynamics through damping and predictive measures.  We apply our PID Lagrangian methods in deep RL, setting a new state of the art in Safety Gym, a safe RL benchmark.  Lastly, we introduce a new method to ease controller tuning by providing invariance to the relative numerical scales of reward and cost.  Our extensive experiments demonstrate improved performance and hyperparameter robustness, while our algorithms remain nearly as simple to derive and implement as the traditional Lagrangian approach.


\end{abstract}

\section{Introduction}
\label{sec:intro}

Reinforcement learning has solved sequential decision tasks of impressive difficulty by maximizing reward functions through trial and error.  Recent examples using deep learning range from robotic locomotion \cite{trpo, ddpg, PPO, levine2016} to sophisticated video games \cite{DQN, PPO, OpenAI_dota, jaderberg2019human}.  While errors during training in these domains come without cost, in some learning scenarios it is important to limit the rates of hazardous outcomes.  One example would be wear and tear on a robot's components or its surroundings.  It may not be possible to impose such limits by prescribing constraints in the action or state space directly; instead, hazard-avoiding behavior must be learned.  For this purpose, we use the well-known framework of the constrained Markov decision process (CMDP) \cite{altman1999}, which limits the accumulation of a ``cost'' signal which is analogous to the reward.  The optimal policy is one which maximizes the usual return while satisfying the cost constraint.  In \emph{safe} RL the agent must avoid hazards not only at convergence, but also throughout exploration and learning.

Lagrangian methods are a classic approach to solving constrained optimization problems.  For example, the equality-constrained problem over the real vector $\mathbf{x}$:
\begin{equation}
    \min_\mathbf{x} f(\mathbf{x}) \;\;\;\;\; \text{s.t.} \; g(\mathbf{x})=0
\end{equation}
is transformed into an unconstrained one by introduction of a dual variable--the Lagrange multiplier, $\lambda$--to form the Lagrangian: $\mathcal{L}(\mathbf{x}, \lambda) = f(\mathbf{x}) + \lambda g(\mathbf{x})$,
which is used to find the solution as:
\begin{equation}
    (\mathbf{x}^*,\lambda^*)=\arg\max_\lambda\min_\mathbf{x} \mathcal{L}(\mathbf{x},\lambda)
\end{equation}
Gradient-based algorithms iteratively update the primal and dual variables:
\begin{align}
  -\nabla_\mathbf{x}\mathcal{L}(\mathbf{x},\lambda)=&-\nabla_\mathbf{x} f(\mathbf{x}) - \lambda \nabla_\mathbf{x} g(\mathbf{x})\\
  \nabla_\lambda \mathcal{L}(\mathbf{x},\lambda)=&\ g(\mathbf{x})
  \label{eq:nabla_lambda}
\end{align}
so that $\lambda$ acts as a learned penalty coefficient in the objective, leading eventually to a constraint-satisfying solution (see \emph{e.g.} \citet{bertsekas2014constrained}).  The Lagrangian multiplier method is readily adapted to the constrained RL setting \cite{altman1998, geibel_rl} and has become a popular baseline in deep RL \cite{CPO, chow2019} for its simplicity and effectiveness.

\begin{figure}[ht]
\begin{center}
\centerline{\includegraphics[width=1.05\columnwidth]{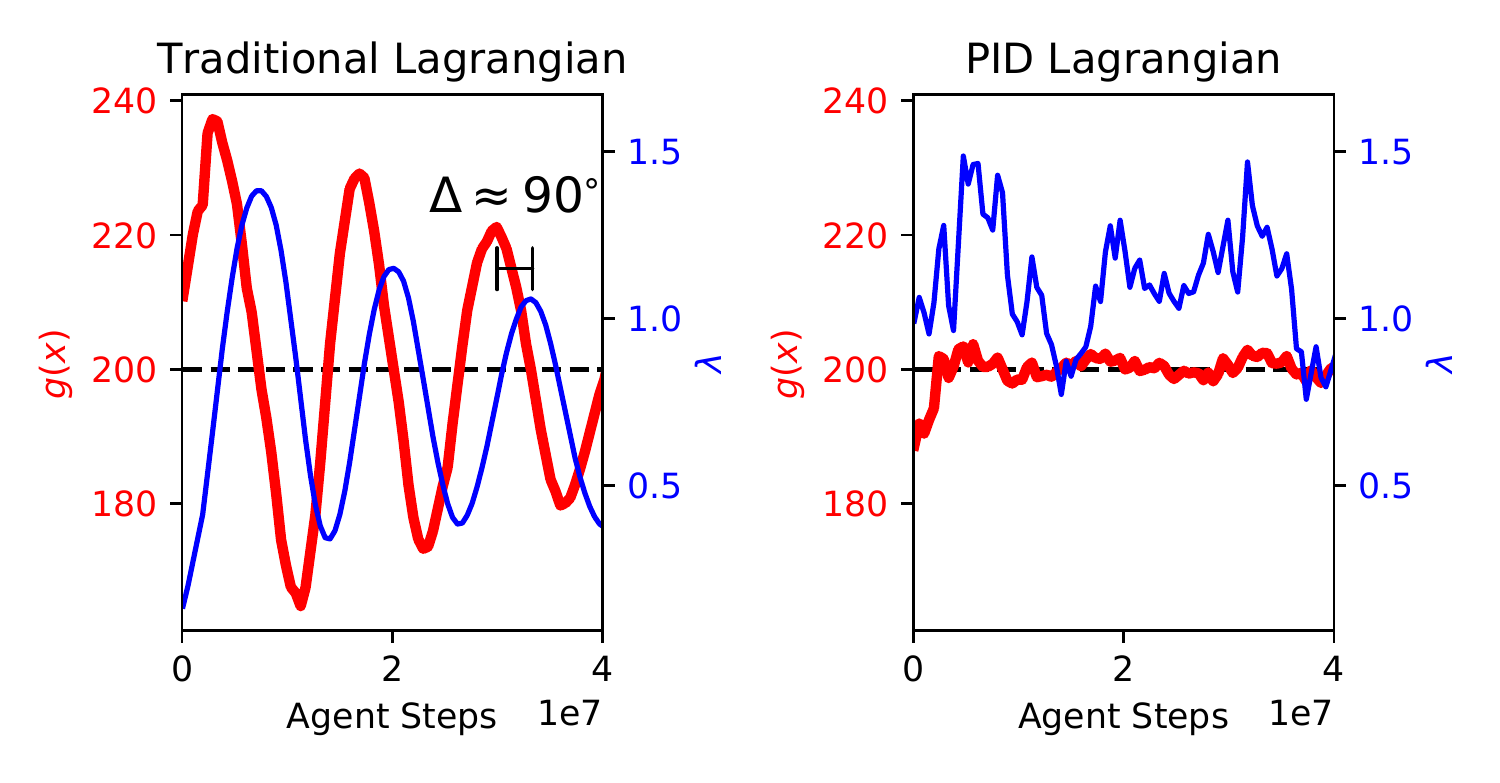}}
\vskip -0.2in
\caption{\emph{Left}: The traditional Lagrangian method exhibits oscillations with 90$^{\circ}$ phase shift between the constraint function and the Lagrange multiplier, characteristic of integral control. \emph{Right}: PID control on the Lagrange multiplier damps oscillations and obeys constraints. Environment: \textsc{DoggoButton1}, cost limit 200.}
\label{fig:phase_delay}
\end{center}
\vskip -0.2in
\end{figure}

Although they have been shown to converge to optimal, constraint-satisfying policies \cite{rcpo, paternain_duality_gap}, a shortcoming of gradient Lagrangian methods for safe RL is that intermediate iterates often violate constraints.  Cost overshoot and oscillations are in fact inherent to the learning dynamics \cite{platt1988constrained, wah2000improving}, and we witnessed numerous problematic cases in our own experiments.  Figure~\ref{fig:phase_delay} (\emph{left}) shows an example from a deep RL setting, where the cost and multiplier values oscillated throughout training.  Our key insight in relation to this deficiency is that the traditional Lagrange multiplier update in (\ref{eq:nabla_lambda}) amounts to \emph{integral} control on the constraint.  The 90-degree phase shift between the curves is characteristic of ill-tuned integral controllers.

Our contribution is to expand the scope of possible Lagrange multiplier update rules beyond (\ref{eq:nabla_lambda}), by interpreting the overall learning algorithm as a dynamical system.  Specifically, we employ the next simplest mechanisms, \emph{proportional} and \emph{derivative} control, to $\lambda$, by adding terms corresponding to derivatives of the constraint function into (\ref{eq:nabla_lambda}) (derivatives with respect to learning iteration).  To our knowledge, this is the first time that an expanded update rule has been considered for a learned Lagrange multiplier.  PID control is an appealing enhancement, evidenced by the fact that it is one of the most widely used and studied control techniques \cite{aastrom2006pid}.  The result is a more responsive safety mechanism, as demonstrated in Figure~\ref{fig:phase_delay} (\emph{right}), where the cost oscillations have been damped, dramatically reducing violations.

Our contributions in this paper are outlined as follows.  First, we provide further context through related works and preliminary definitions.  In Section~\ref{sec:lagrange}, we propose modified Lagrangian multiplier methods and analyze their benefits in the learning dynamics.  Next, in Section~\ref{sec:feedback}, we cast constrained RL as a dynamical system with the Lagrange multiplier as a control input, to which we apply PID control as a new algorithm.  In Section~\ref{sec:experiments}, we adapt a leading deep RL algorithm, Proximal Policy Optimization (PPO) \cite{PPO} with our methods and achieve state of the art performance in the OpenAI Safety-Gym suite of environments \cite{safetygym}.  Finally, in Section~\ref{sec:invariance} we introduce another novel technique that makes tuning easier by providing invariance to the relative numerical scales of rewards and costs, and we demonstrate it in a further set of experiments.  Our extensive empirical results show that our algorithms, which are intuitive and simple to implement, improve cost performance and promote hyperparameter robustness in a deep RL setting.

\section{Related Work}
\label{sec:related}
\textbf{Constrained Deep RL.} Adaptations of the Lagrange multiplier method to the actor-critic RL setting have been shown to converge to the optimal, constraint-satisfying solution under certain assumptions \cite{rcpo}.  Convergence proofs have relied upon updating the multiplier more slowly than the policy parameters \cite{rcpo, paternain_duality_gap}, implying many constraint-violating policy iterations may occur before the penalty comes into full effect.

Several recent works have aimed at improving constraint satisfaction in RL over the Lagrangian method, but they tend to incur added complexity.  \citet{CPO} introduced Constrained Policy Optimization (CPO), a policy search algorithm with near-constraint satisfaction guarantees at every iteration, based on a new bound on the expected returns of two nearby policies.  CPO includes a projection step on the policy parameters, which in practice requires a time-consuming backtracking line search.  Yet, simple Lagrangian-based algorithms performed as well or better in a recent empirical comparison in Safety Gym \cite{safetygym}. Approaches to safe RL based on Lyapunov functions have been developed in a series of studies \cite{Chow2018, chow2019}, resulting in algorithms that combine a projection step, as in CPO, with action-layer interventions like the safety layer of \citet{safety-layer}. Experimentally, this line of work showed mixed performance gains over Lagrangian methods, at a nontrivial cost to implement and without clear guidance for tuning. \citet{liu2019ipo} developed interior point methods for RL, which augment the objective with logarithmic barrier functions.  These methods are shown theoretically to provide suboptimal solutions.  Furthermore, they require tuning of the barrier strength and typically assume already feasible iterates, the latter point possibly being problematic for random agent initializations or under noisy cost estimates.  Most recently, \citet{Yang2020Projection-Based} extended CPO with a two-step projection-based optimization approach.  In contrast to these techniques, our method remains nearly as simple to implement and compute as the baseline Lagrangian method.

\textbf{Dynamical Systems View of Optimization.}   
Several recent works have proposed different dynamical systems viewpoints to analyze optimization algorithms, including those often applied to deep learning. \citet{hu_lessard} reinterpreted first-order gradient optimization as a dynamical system; they likened the gradient of the objective, $\nabla_x f$, to the plant, which the controller aims to drive to zero to arrive at the optimal parameters, $x^*$.  Basic gradient descent then matches the form of integral control (on $\nabla_x f$).  They extend the analogy to momentum-based methods, for example linking Nesterov momentum to PID control with lag compensation.  In another example, \citet{An2018PID} interpreted SGD as P-control and momentum methods as PI-control.  They introduced a derivative term, based on the change in the gradient, and applied their resulting PID controller to improve optimization of deep convolutional networks.  Other recent works bring yet other perspectives from dynamical systems to deep learning and optimization, see for example \cite{lessard2014analysis, nishihara2015general, liu2019deep}).  None of these works address constrained RL, however, necessitating our distinct formulation for that problem.

\textbf{Constrained Optimization.}
Decades' worth of literature have accumulated on Lagrangian methods.  But even recent textbooks on the topic \cite{bertsekas2014constrained,nocedal2006numerical} only consider updating the Lagrange multiplier using the value of the constraint function, $g(x)$, and miss ever using its derivatives, $\dot{g}(x)$ or $\ddot{g}(x)$, which we introduce.  The modification to the Lagrangian method most similar in effect to our proportional control term (here using $\dot{g}(x)$) is the quadratic penalty method (\citet{hestenes1969multiplier, powell1969method} see also \textit{e.g.} \citet{bertsekas1976penalty}), which we compare in Section~\ref{sec:lagrange}.  \citet{Song1998} proposed a controls viewpoint (continuous-time) of optimizing neural networks for constrained problems and arrived at proportional control rules only.  Related to our final experiments on reward-scale invariance, \citet{wah2000improving} developed an adaptive weighting scheme for continuous-time Lagrangian objectives, but it is an intricate procedure which is not straightforwardly applied to safe RL.

\section{Preliminaries}
\label{sec:prelim}

\textbf{Constrained Reinforcement Learning}  Constrained Markov Decision Processes (CMDP) \cite{altman1998} extend MDPs (see \citet{suttonbarto}) to incorporate constraints into reinforcement learning.  A CMDP is the expanded tuple $(S,A,R,T,\mu,C_0,C_1,...,d_0,d_1,...)$, with the cost functions $C_i:S\times A\times S \rightarrow \mathbb{R}$ defined by the same form as the reward, and $d_i: \mathbb{R}$ denoting limits on the costs.  For ease of notation, we will only consider a single, all-encompassing cost.  

The expected sum of discounted rewards over trajectories, $\tau=(s_0,a_0,s_1,a_1,...)$, induced by the policy $\pi(a|s)$ is a common performance objective: $J(\pi) = \mathbb{E}_{\tau\sim\pi}\left[\sum_{t=0}^{\infty}{\gamma^t R(s_t,a_t,s_{t+1})}\right]$.  The analogous value function for cost is defined as: $J_{C}(\pi)=\mathbb{E}_{\tau\sim\pi}\left[\sum_{t=0}^\infty{\gamma^t C(s_t,a_t,s_{t+1})}\right]$.  The constrained RL problem is to solve for the best feasible policy:
\begin{equation}
\label{eq:constrained_rl}
     \pi^* = \arg \max_{\pi} J(\pi) \;\;\;\;\; \text{s.t.} \; J_{C}(\pi) \leq d
\end{equation}

Deep reinforcement learning uses a deep neural network for the policy, $\pi_\theta=\pi(\cdot|s;\theta)$ with parameter vector $\theta$, and policy gradient algorithms improve the policy iteratively by gathering experience in the task to estimate the reward objective gradient, $\nabla_\theta J(\pi_\theta)$.  Thus our problem of interest is better expressed as maximizing score at some iterate, $\pi_k$, while ideally obeying constraints at each iteration:
\begin{equation}
\label{eq:P1}
    \begin{aligned}
    &\max_\pi J(\pi_k)\\
    &\text{s.t.}\; J_C(\pi_m)\leq d \quad \; m\in\{0,1,...,k\}
    \end{aligned}
\end{equation}
Practical settings often allow trading reward performance against some constraint violations (\emph{e.g.} the constraints themselves may include a safety margin).  For this purpose we introduce a constraint figure of merit with our experiments.

\subsection{Dynamical Systems and Optimal Control}
Dynamical systems are processes which can be subject to an external influence, or \textit{control}.  A general formulation for discrete-time systems with feedback control is:
\begin{equation}
\label{eq:dynamical}
\begin{aligned}
    \mathbf{x}_{k+1} =& F(\mathbf{x}_k,\mathbf{u}_k)\\
    \mathbf{y}_k =& Z(\mathbf{x}_k)\\
    \mathbf{u}_k =& h(\mathbf{y}_0,...,\mathbf{y}_k)
\end{aligned}
\end{equation}
with state vector $\mathbf{x}$, dynamics function $F$, measurement outputs $\mathbf{y}$, applied control $\mathbf{u}$, and the subscript denoting the time step.  The feedback rule $h$ has access to past and present measurements.  A problem in optimal control is to design a control rule, $h$, that results in a sequence $\mathbf{y}_{0:T}\doteq \{\mathbf{y}_0,...,\mathbf{y}_T\}$ (or $\mathbf{x}_{0:T}$ directly) that scores well according to some cost function $C$. Examples include simply reaching a goal condition, $C=|\mathbf{y}_T-\overline{\mathbf{y}}|$, or following close to a desired trajectory, $\overline{\mathbf{y}}_{0:T}$.
   
Systems with simpler dependence on the input are generally easier to analyze and control (\emph{i.e.} simpler $h$ performs well), even if the dependence on the state is complicated \cite{skelton1988dynamic}.  Control-affine systems are a broad class of dynamical systems which are especially amenable to analysis \cite{isidori}.  They take the form:
\begin{equation}
\label{eq:affine}
    F(\mathbf{x}_k,\mathbf{u}_k)=f(\mathbf{x}_k) + g(\mathbf{x}_k)\mathbf{u}_k
\end{equation}
where $f$ and $g$ may be nonlinear in state, and are possibly \textit{uncertain}, meaning unknown.  We will seek control-affine form for ease of control and to support future analysis.  

\section{Modified Lagrangian Methods for Constrained Optimization}
\label{sec:lagrange}

Lagrangian methods are a classic family of approaches to solving constrained optimization problems.  We propose an intuitive, previously overlooked form for the multiplier update and derive its beneficial effect on the learning dynamics.  We begin by reviewing a prior formulation for the equality-constrained problem.\footnote{Standard techniques extend our results to inequality constraints, and multiple constraints, as in \citet{platt1988constrained}, and notation is simplest for an equality constraint.}


\subsection{Review: ``Basic Differential Multiplier Method''}
We follow the development of \citet{platt1988constrained}, who analyzed the dynamics of a continuous-time neural learning system applied to this problem (our result can similarly be derived for iterative gradient methods).  They begin with the component-wise differential equations:
\begin{align}
 \dot{x_i}&=-\frac{\partial{\mathcal{L}(\mathbf{x},\lambda)}}{\partial{x_i}}=-\frac{\partial{f}}{\partial{x_i}}-\lambda\frac{\partial{g}}{\partial{x_i}} \label{eq:x_i_dot} \\ 
 \Dot{\lambda}&=\alpha \frac{\partial{\mathcal{L}(\mathbf{x},\lambda)}}{\partial{\lambda}}=\alpha g(\mathbf{x}) \label{eq:lambda_dot}
\end{align}
where we have inserted the scalar constant $\alpha$ as a learning rate on $\lambda$.  Differentiating (\ref{eq:x_i_dot}) and substituting with (\ref{eq:lambda_dot}) leads to the second-order dynamics, written in vector format:
\begin{equation}
\label{eq:I_dynamics}
    \ddot{\mathbf{x}}+A\dot{\mathbf{x}}+\alpha g(\mathbf{x})\nabla g=0
\end{equation}
which is a forced oscillator with damping matrix equal to the weighted sum of Hessians:
\begin{equation}
    A_{ij}=\frac{\partial^2 f}{\partial x_i \partial x_j} + \lambda \frac{\partial^2 g}{\partial x_i \partial x_j}, \text{ or, } A = \nabla^2 f + \lambda \nabla^2 g
\end{equation}
\citet{platt1988constrained} showed that if $A$ is positive definite, then the system (\ref{eq:I_dynamics}) converges to a solution that satisfies the constraint. \citet{platt1988constrained} also noted that the system (\ref{eq:x_i_dot})-(\ref{eq:lambda_dot}) is prone to oscillations as it converges into the feasible region, with frequency and settling time depending on $\alpha$.  We provide complete derivations of the dynamics in (\ref{eq:I_dynamics}) and for our upcoming methods in an appendix. 
  
\subsection{Proportional-Integral Multiplier Method}
In (\ref{eq:lambda_dot}), $\lambda$ simply integrates the constraint.  To improve the dynamics towards more rapid and stable satisfaction of constraints, we introduce a new term in $\lambda$ that is \emph{proportional} to the current constraint value.  In the differential equation for $\lambda$, this term appears as the time-derivative of the constraint:
\begin{equation}
\label{eq:PI_lambda_dot}
    \Dot{\lambda}= \alpha g(\mathbf{x}) + \beta \dot{g}(\mathbf{x})= \alpha g(\mathbf{x}) + \beta \sum_j{\frac{\partial g}{\partial x_j}}\dot{x}_j
\end{equation}
with strength coefficient, $\beta$.  Replacing (\ref{eq:lambda_dot}) by (\ref{eq:PI_lambda_dot}) and combining with (\ref{eq:x_i_dot}) yields similar second-order dynamics as (\ref{eq:I_dynamics}), except with an additional term in the damping matrix:
\begin{equation}
\label{eq:PI_dynamics}
    \ddot{\mathbf{x}}+\left(A + \beta \nabla g \nabla^\top g \right)\dot{\mathbf{x}}+\alpha g(\mathbf{x})\nabla g=0
\end{equation}
The new term is beneficial because it is positive semi-definite---being the outer product of a vector with itself---so it can increase the damping eignevalues, boosting convergence.  The results of \cite{platt1988constrained} hold under (\ref{eq:PI_lambda_dot}, \ref{eq:PI_dynamics}), because the conditions of the solution, namely $\dot{\mathbf{x}}=0$ and $g(\mathbf{x})=0$, remain unaffected and extend immediately to $\dot{g}(\mathbf{x})=0$ (and for the sequel, to $\ddot{g}(\mathbf{x})=0$).  To our knowledge, this is the first time that a proportional-integral update rule has been considered for a learned Lagrange multiplier.  

The well-known penalty method \cite{hestenes1969multiplier,powell1969method} augments the Lagrangian with an additional term, $\frac{c}{2}g(\mathbf{x})^2$, which produces a similar effect on the damping matrix, as shown in \cite{platt1988constrained}:
\begin{equation}
    A_{penalty} = A + c \nabla g \nabla^\top g + c g(\mathbf{x})\nabla^2 g
\end{equation}
Our approach appears to provide the same benefit, without the following two complications of the penalty method.  First, the penalty term must be implemented in the derivative $\dot{\mathbf{x}}$, whereas our methods do not modify the Lagrangian nor the derivative in (\ref{eq:x_i_dot}).  Second, the penalty introduces another instance of the hessian $\nabla^2 g$ in the damping matrix, which might not be positive semi-definite but shares the proportionality factor, $c$, with the desired term.


\subsection{Integral-Derivative Multiplier Method}
A similar analysis extends to the addition of a term in $\lambda$ based on the \emph{derivative} of the constraint value.  It appears in $\dot{\lambda}$ as the second derivative of the constraint:
\begin{equation}
    \dot{\lambda} = \alpha g(\mathbf{x}) + \gamma \ddot{g}(\mathbf{x})
\end{equation}
with strength coefficient $\gamma$. The resulting dynamics are:
\begin{equation}
    \label{eq:ID_dynamics}
    \ddot{\mathbf{x}} + B^{-1}A\dot{\mathbf{x}} + \left(\alpha g(\mathbf{x}) + \gamma \dot{\mathbf{x}}^\top \nabla^2 g \dot{\mathbf{x}}  \right)  B^{-1} \nabla g = 0
\end{equation}
with $B=\left(I + \gamma \nabla g \nabla^\top g\right)$, and $I$ the identity matrix.  

The effects of the derivative update method are two-fold.  First, since the eigenvalues of the matrix $B^{-1}$ will be less than $1$, both the damping ($A$) and forcing ($\nabla g$) terms are weakened (and rotated, generally).  Second, the new forcing term can be interpreted as a drag quadratic in the speed and modulated by the curvature of the constraint along the direction of motion.  To illustrate cases, if the curvature of $g$ is positive along the direction of travel, then this term becomes a force for decreasing $g$.  If at the same time $g(\mathbf{x})>0$, then the traditional force will also be directed to decrease $g$, so the two will add.  On the other hand, if $g$ curves negatively along the velocity, then the new force promotes increasing $g$; if $g(\mathbf{x})>0$, then the two forces subtract, weakening the acceleration $\ddot{\mathbf{x}}$.  By using curvature, the derivative method acts predictively, but may be prone to instability.

The \textbf{proportional-integral-derivative multiplier method} is the combination of the previous two developments, which induced independent changes in the dynamics (\emph{i.e.} insert the damping matrix of (\ref{eq:PI_dynamics}) into (\ref{eq:ID_dynamics})).  We leave for future work a more rigorous analysis of the effects of the new terms, along with theoretical considerations of the values of coefficients $\alpha$, $\beta$, and $\gamma$.  In the next section, we carry the intuitions from our analysis to make practical enhancements to Lagrangian-based constrained RL algorithms.

\section{Feedback Control for Constrained RL}
\label{sec:feedback}
We advance the broader consideration of possible multiplier update rules by reinterpreting constrained RL as a dynamical system; the adaptive penalty coefficient is a control input, and the cost threshold is a setpoint which the system should maintain.  As the agent learns for rewards, the upward pressure on costs from reward-learning can change, requiring dynamic response.  In practical Lagrangian RL, the iterates $\lambda_k$ may deviate from the optimal value, even for lucky initialization $\lambda_0=\lambda^*$, as the policy is only partially optimized at each iteration.  Adaptive sequences $\lambda_0,...,\lambda_K$ other than those prescribed by the Lagrangian method may achieve superior cost control for Problem (\ref{eq:P1}).  In this section we relate the Lagrangian method to a dynamical system, formalizing how to incorporate generic update rules using feedback.  We return to the case of an inequality constrained CMDP to present our main algorithmic contribution---the use of PID control to adapt the penalty coefficient.

\subsection{Constrained RL as a Dynamical System}
We write constrained RL as the first-order dynamical system:
\begin{equation}
    \begin{aligned}
        \theta_{k+1}=&F(\theta_k, \lambda_k)\\
        y_k=&J_C(\pi_{\theta_k})\\
        \lambda_k=&h(y_0,...,y_k, d)
    \end{aligned}
    \label{eq:constrained_rl_dynamical}
\end{equation}
where $F$ is an unknown nonlinear function\footnote{Known as an ``uncertain'' nonlinear function in the control literature, meaning we lack an analytical expression for it.} corresponding to the RL algorithm policy update on the agent's parameter vector, $\theta$.  The cost-objective serves as the system measure, $y$, which is supplied to the feedback control rule, $h$, along with cost limit, $d$.  From this general starting point, both the RL algorithm, $F$, and penalty coefficient update rule, $h$, can be tailored for solving Problem (\ref{eq:P1}).

The reward and cost policy gradients of the first-order\footnote{We discuss only the first-order case, which provides sufficient clarity for our developments.} Lagrangian method, $\nabla_\theta\mathcal{L}(\theta,\lambda)=\nabla_\theta J(\pi_\theta)-\lambda \nabla_\theta J_C(\pi_\theta)$, can be organized into the form of (\ref{eq:constrained_rl_dynamical}) as:
\begin{equation}
\label{eq:F_affine}
    F(\theta_k,\lambda_k)= f(\theta_k)+g(\theta_k)\lambda_k
\end{equation}
\begin{equation}
\label{eq:f_GD}
    f(\theta_k)= \theta_k + \eta\nabla_\theta J(\pi_{\theta_k})
\end{equation}
\begin{equation}
\label{eq:g_lagrange}
    g(\theta_k)=-\eta\nabla_\theta J_C(\pi_{\theta_k})
\end{equation}
with SGD learning rate $\eta$.  The role of the controller is to drive inequality constraint violations $(J_c-d)_+$ to zero in the presence of drift from reward-learning due to $f$.  The Lagrange multiplier update rule for an inequality constraint uses subgradient descent:
\begin{equation}
\label{eq:lambda_old}
    \lambda_{k+1}=(\lambda_k + K_I(J_C-d))_+    
\end{equation}
with learning rate $K_I$ and projection into $\lambda\geq 0$.  This update step is clearly an \emph{integral} control rule, for $h$.

\subsection{Constraint-Controlled RL}
Our general procedure, constraint-controlled RL, is given in Algorithm \ref{algo:constrained_rl}.  It follows the typical minibatch-RL scheme, and sampled estimates of the cost criterion, $\hat{J}_C$ are fed back to control the Lagrange multiplier.  In contrast to prior work \cite{rcpo, paternain_duality_gap} which uses a single value approximator and treats $r+\lambda c$ as the reward, we use separate value- and cost-value approximators, since $\lambda$ may change rapidly.

When $\lambda$ is large, the update in (\ref{eq:F_affine}) can cause excessively large change in parameters, $\theta$, destabilizing learning.  To maintain consistent step size, we use a re-scaled objective for the $\theta$-learning loop:
\begin{equation*}
    \theta^*(\lambda)=\arg\max_\theta J-\lambda J_C=\arg\max_\theta \frac{1}{1+\lambda} (J-\lambda J_C)
\end{equation*}
This convex combination of objectives yields the policy gradient used in Algorithm \ref{algo:constrained_rl}.  Our experiments use this re-scaling, including for traditional Lagrangian baselines. 

\begin{algorithm}
\caption{\\Constraint-Controlled Reinforcement Learning}\label{algo:constrained_rl}
\begin{algorithmic}[1]

\Procedure{Constrained RL}{$\pi_{\theta_0}(\cdot|s), d$}
\State Initialize control rule (as needed)
\State $\mathcal{J_C} \gets \{\}$  \Comment{cost measurement history}
\Repeat

\State Sample environment: \Comment{a minibatch}
\State $\quad a\sim\pi(\cdot|s;\theta), s'\sim T(s,a),$
\State $\quad r\sim R(s,a,s'),c\sim C(s,a,s')$

\State Apply feedback control:
\State $\quad$Store sample estimate $\hat{J_C}$ into $\mathcal{J_C}$
\State $\quad \lambda\gets h(\mathcal{J_C}, d),\;\lambda \geq 0$

\State Update $\pi$ by RL: \Comment{by Lagrangian objective}
\State $\quad$ Update critics, $V_\phi(s), V_{C,\psi}(s)$  \Comment{if using}
\State $\quad$ $\nabla_\theta \mathcal{L}=\frac{1}{1+\lambda}\left(\nabla_\theta\hat{J}(\pi_\theta) -\lambda \nabla_\theta\hat{J_C}(\pi_\theta)\right)$

\Until{converged}
\State \textbf{return} $\pi_\theta$
\EndProcedure
\end{algorithmic}
\end{algorithm}

As an aside, we note that it is possible to maintain the control-affine form of (\ref{eq:F_affine}) with this re-scaling, by reparameterizing the control as $0\leq u=\frac{\lambda}{1+\lambda}\leq 1$ and substituting for (\ref{eq:g_lagrange}) with:
\begin{equation}
\label{eq:g_CC}
    g(\theta_k)=-\eta\nabla_\theta\left(J(\pi_{\theta_k}) + J_C(\pi_{\theta_k})\right)
\end{equation}
This parameterization simply weights the reward and cost gradients in the Lagrangian objective as:
\begin{equation}
\nabla_\theta\mathcal{L}(\theta,\lambda)=(1-u)\nabla_\theta J(\pi_\theta)-u \nabla_\theta J_C(\pi_\theta)    
\end{equation}
It may provide superior performance in some cases, as it will behave differently in relation to the nonlinearity in control which arises from the inequality constraint.  We leave experimentation with direct control on $u\in[0,1]$ to future work.

\subsection{The PID Lagrangian Method}
We now specify a new control rule for use in Algorithm \ref{algo:constrained_rl}.  
To overcome the shortcomings of integral-only control, we follow the developments of the previous section and introduce the next simplest components: \emph{proportional} and \emph{derivative} terms.  Our PID update rule to replace (\ref{eq:lambda_old}) is shown in Algorithm \ref{algo:lagrange_pid}.  The proportional term will hasten the response to constraint violations and dampen oscillations, as derived in Section \ref{sec:lagrange}.  Unlike the Lagrangian update, derivative control can act in anticipation of violations.  It can both prevent cost overshoot and limit the rate of cost increases within the feasible region, useful when monitoring a system for further safety interventions.  Our derivative term is projected as $(\cdot)_+$ so that it acts against increases in cost but does not impede decreases.  Overall, PID control provides a much richer set of controllers while remaining nearly as simple to implement; setting $K_P=K_D=0$ recovers the traditional Lagrangian method.  The integral term remains necessary for eliminating steady-state violations at convergence.  Our experiments mainly focus on the effects of proportional and derivative control of the Lagrange multiplier in constrained deep RL.

\begin{algorithm}
\caption{PID-Controlled Lagrange Multiplier}\label{algo:lagrange_pid}
\begin{algorithmic}[1]
\State Choose tuning parameters: $K_P, K_I, K_D \geq 0$
\State Integral: $I\gets0$
\State Previous Cost: $J_{C,prev}\gets 0$
\MRepeat \ at each iteration $k$
\State Receive cost $J_C$
\State $\Delta \gets J_C-d$
\State $\partial \gets (J_C-J_{C,prev})_+$
\State $I \gets (I + \Delta)_+$
\State $\lambda \gets (K_P \Delta + K_I I + K_D \partial)_+$
\State $J_{C,prev}\gets J_C$
\State \textbf{return} $\lambda$
\EndRepeat
\end{algorithmic}
\end{algorithm}

\section{PID Control Experiments}
\label{sec:experiments}

We investigated the performance of our algorithms on Problem (\ref{eq:P1}) in a deep RL setting.  In particular, we show the effectiveness of PID control at reducing constraint violations from oscillations and overshoot present in the baseline Lagrangian method.  Both maximum performance and robustness to hyperparameter selection are considered.  Although many methods exist for tuning PID parameters, we elected to do so manually, demonstrating ease of use. 


\subsection{Environments: Safety-Gym}
We use the recent Safety-Gym suite \cite{safetygym}, which consists of robot locomotion tasks built on the MuJoCo simulator \cite{mujoco}.  The robots range in complexity from a simple Point robot to the 12-jointed Doggo, and they move in an open arena floor.  Rewards have a small, dense component encouraging movement toward the goal, and a large, sparse component for achieving it.  When a goal is achieved, a new goal location is randomly generated, and the episode continues until the time limit at 1,000 steps.

\begin{figure}[ht]
\begin{center}
\centerline{
  \includegraphics[clip,width=0.6\columnwidth]{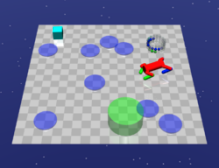}}
\vskip -0.1in
\caption{Rendering from the \textsc{DoggoGoal1} environment from Safety Gym.  The red, four-legged robot must walk to the green cylinder while avoiding other objects, and receives coarse egocentric sensor readings of their locations.}
\label{fig:doggogoal_screenshot}
\end{center}
\vskip -0.2in
\end{figure}

Each task has multiple difficulty levels corresponding to density and type of hazards, which induce a cost when contacted by the robot (without necessarily hindering its movement).  Hazards are placed randomly at each episode and often lay in the path to the goal.  Hence the aims of achieving high rewards and low costs are in opposition.  The robot senses the position of hazards and the goal through a coarse, LIDAR-like mode.  The output of this sensor, along with internal readings like the joint positions and velocities, comprises the state fed to the agent.  Figure~\ref{fig:doggogoal_screenshot} displays a scene from the \textsc{DoggoGoal1} environment.

\subsection{Algorithm: Constraint-Controlled PPO}

We implemented Algorithm \ref{algo:constrained_rl} on top of Proximal Policy Optimization (PPO) \cite{PPO} to make constraint-controlled PPO (CPPO).  CPPO uses an analogous clipped surrogate objective for the cost as for the reward.  Our policy is a 2-layer MLP followed by an LSTM with a skip connection.  We applied smoothing to proportional and derivative controls to accommodate noisy estimates.  The environments' finite horizons allowed use of non-discounted episodic costs as the constraint and input to the controller.  Additional training details can be found in supplementary materials, and our implementation is available at \url{https://github.com/astooke/rlpyt/rlpyt/projects/safe}.

\subsection{Main Results}
We compare PID controller performance against the Lagrangian baseline under a wide range of settings.  Plots showing the performance of the unconstrained analogue confirm that constraints are not trivially satisfied, and they appear in supplementary material.

\subsubsection{Robust Safety with PI Control}
\begin{figure}[ht]
\begin{center}
\centerline{
  \includegraphics[clip,width=1.\columnwidth]{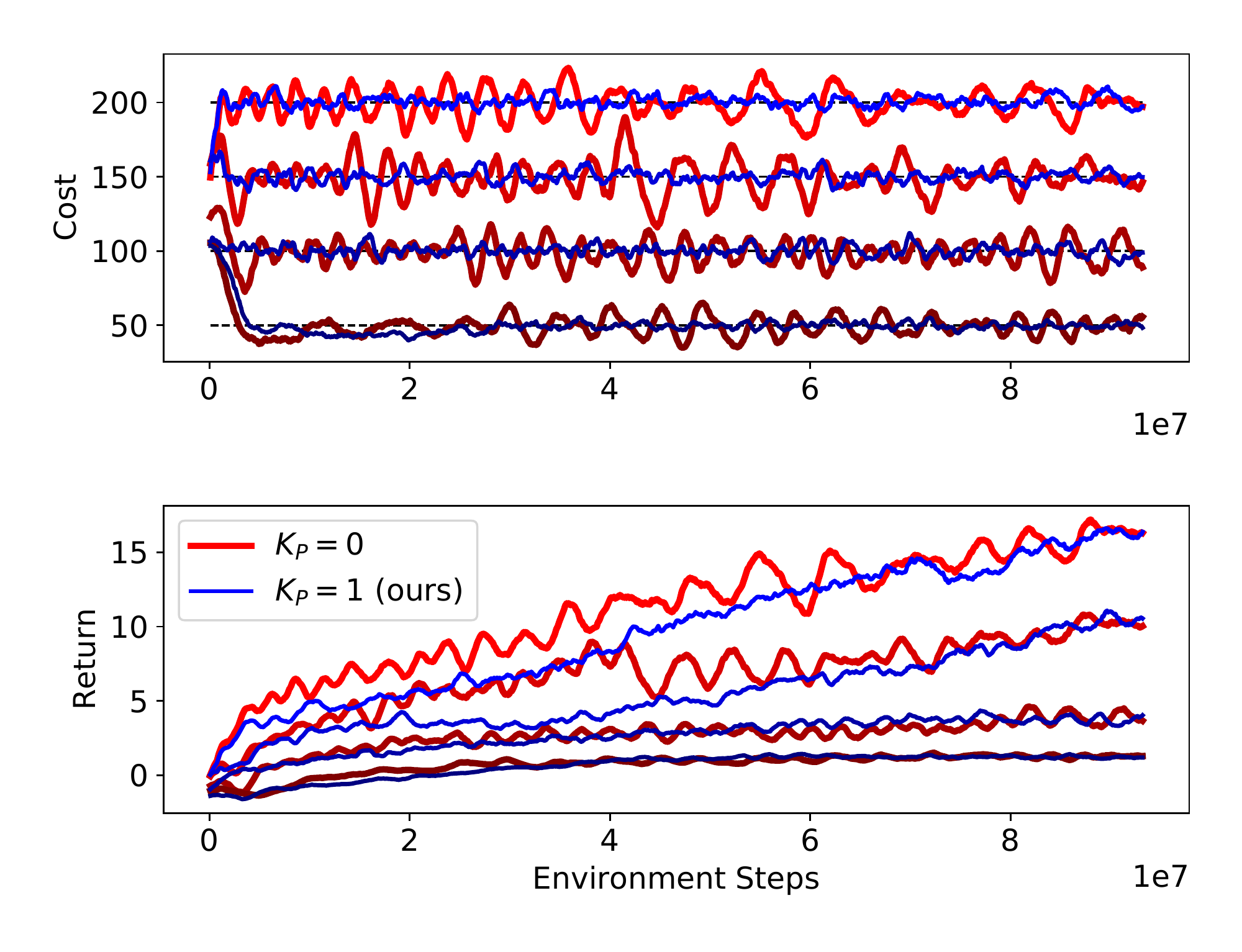}}
\vskip -0.1in
\caption{Oscillations in episodic costs (and returns) from the Lagrangian method, $K_P=0, K_I=10^{-2}$, are damped by proportional control, $K_P=1$ (ours), at cost limits $50, 100, 150, 200$ (curves shaded) in \textsc{DoggoButton1}.}
\label{fig:PI_dampening}
\end{center}
\vskip -0.1in
\end{figure}

We observed cost oscillations or overshoot with slow settling time in a majority of Safety Gym environments when using the Lagrangian method.  Figure \ref{fig:PI_dampening} shows an example where PI-control eliminated this behavior while maintaining good reward performance, in the challenging \textsc{DoggoButton1} environment.  Individual runs are plotted for different cost limits.  

\begin{figure}[ht]
\begin{center}
\centerline{
  \includegraphics[clip,width=1.\columnwidth]{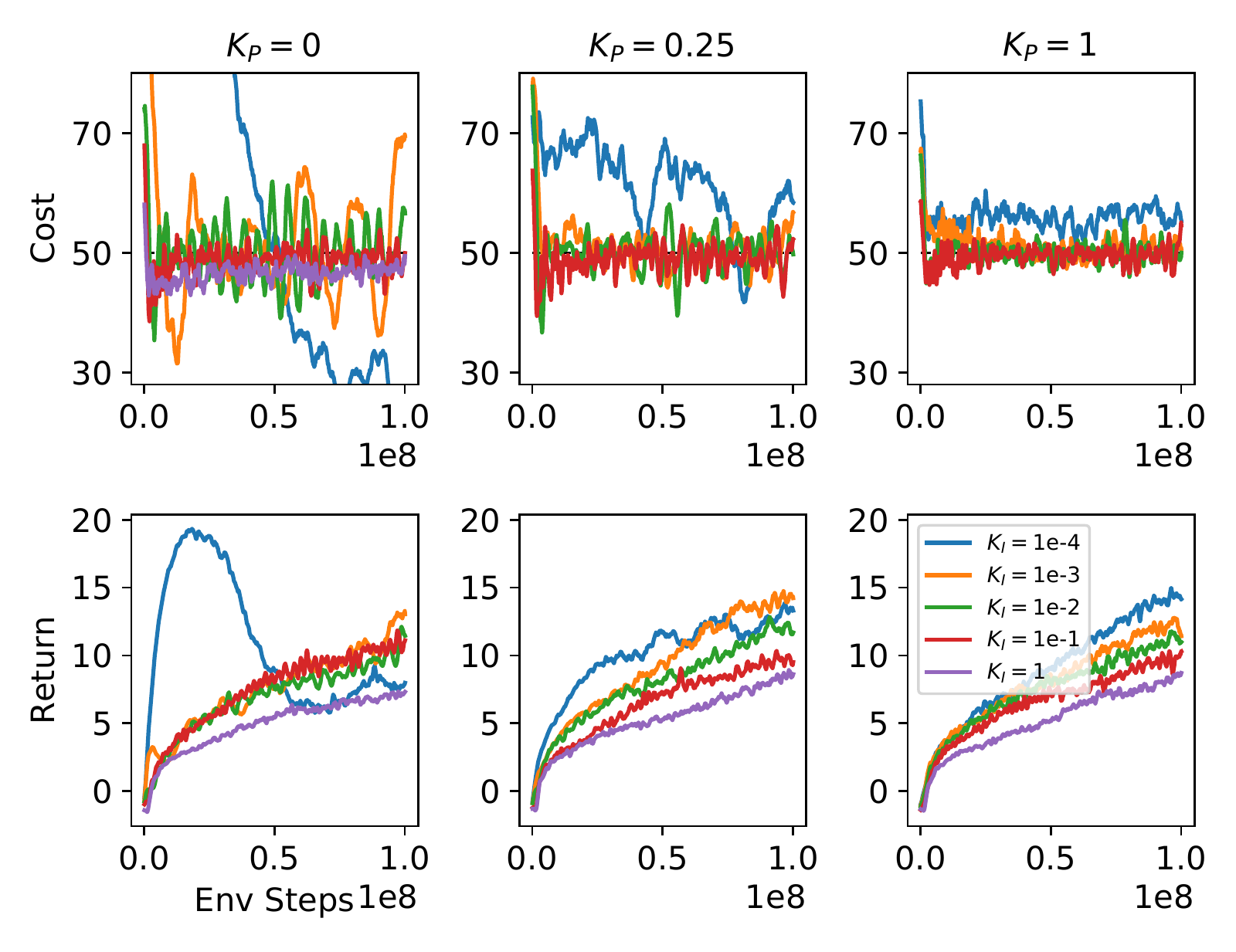}}
\vskip -0.1in
\caption{\textit{Top row}: Constraint-violating oscillations decrease in magnitude and period from increases in the Lagrange multiplier learning rate, $K_I$. At all levels, oscillations are damped by PI-control, $K_P=0.25, 1$.  \textit{Bottom row}: Returns diminish for large $K_I$; proportional control maintains high returns while reducing constraint violations.  Environment: \textsc{DoggoGoal2}, cost limit 50.}
\label{fig:doggogoal_cost_return}
\end{center}
\vskip -0.1in
\end{figure}

As predicted in \cite{platt1988constrained}, we found the severity of cost overshoot and oscillations to depend on the penalty coefficient learning rate, $K_I$.  The top left panel of Figure \ref{fig:doggogoal_cost_return} shows example cost curves from \textsc{DoggoGoal2} under I-control, over a wide range of values for $K_I$ (we refer to varying $K_I$, assuming $K_I=1$; the two are interchangeable in our design).  With increasing $K_I$, the period and amplitude of cost oscillations decrease and eventually disappear.  The bottom left of Figure \ref{fig:doggogoal_cost_return}, however, shows that larger $K_I$ also brings diminishing returns.  We study this effect in the next section.  The center and right columns of Figure \ref{fig:doggogoal_cost_return} show the cost and return when using PI-control, with $K_P=0.25$ and $K_P=1$, respectively.  Proportional control stabilized the cost, with most oscillations reduced to the noise floor for $K_I>10^{-4}$.  Yet returns remained relatively high over a wide range, $K_I<10^{-1}$.  Similar curves for other Safety Gym environments are included in an appendix.

\begin{figure}[ht]

\begin{center}
\centerline{
  \includegraphics[clip,width=0.5\columnwidth]{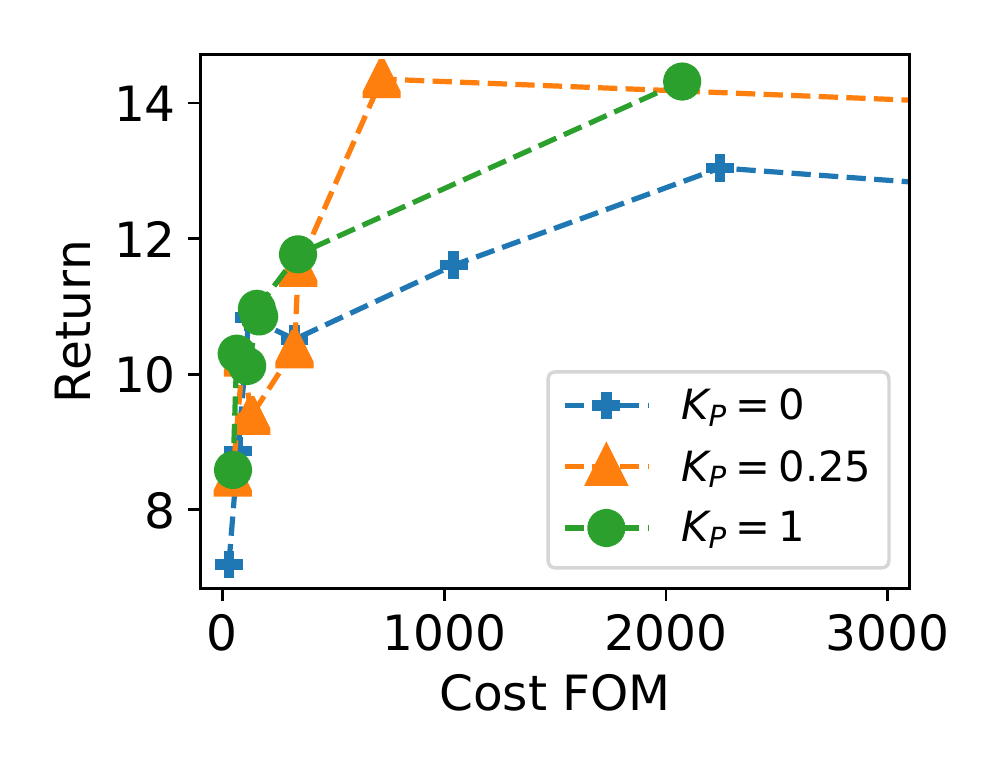}}
\vskip -0.2in
\caption{Pareto frontier of return versus cost FOM, which improves (up and to the left) with PI-control, $K_P=0.25,1$.  Each point is a different setting of $K_I$ (see Figure \ref{fig:doggogoal_cost_return}).}
\label{fig:doggogoal_pareto}
\end{center}
\vskip -0.2in
\end{figure}

We examine the trade-off between reward and constraint violation by forming an overall cost figure of merit (FOM).  We use the sum of non-discounted constraint violations over the learning iterates, $C_{FOM}=\sum_k(D(\pi_{\theta_k})-d)_+, D(\pi_\theta)=\mathbb{E}_{\tau\sim\pi}\left[\sum_{t=0}^T C(s_t,a_t,s_t')\right]$, and estimate it online from the learning data.   Figure \ref{fig:doggogoal_pareto} compares final returns against this cost FOM for the same set of experiments as in Figure \ref{fig:doggogoal_cost_return}.  Each point represents a different setting of $K_I$, averaged over four runs.  PI-control expanded the Pareto frontier of this trade-off into a new region of high rewards at relatively low cost which was inaccessible using the Lagrangian method.  These results constitute a new state of the art over the benchmarks in \citet{safetygym}.

\begin{figure}[ht]
\vskip -0.1in
\begin{center}
\centerline{
  \includegraphics[clip,width=0.9\columnwidth]{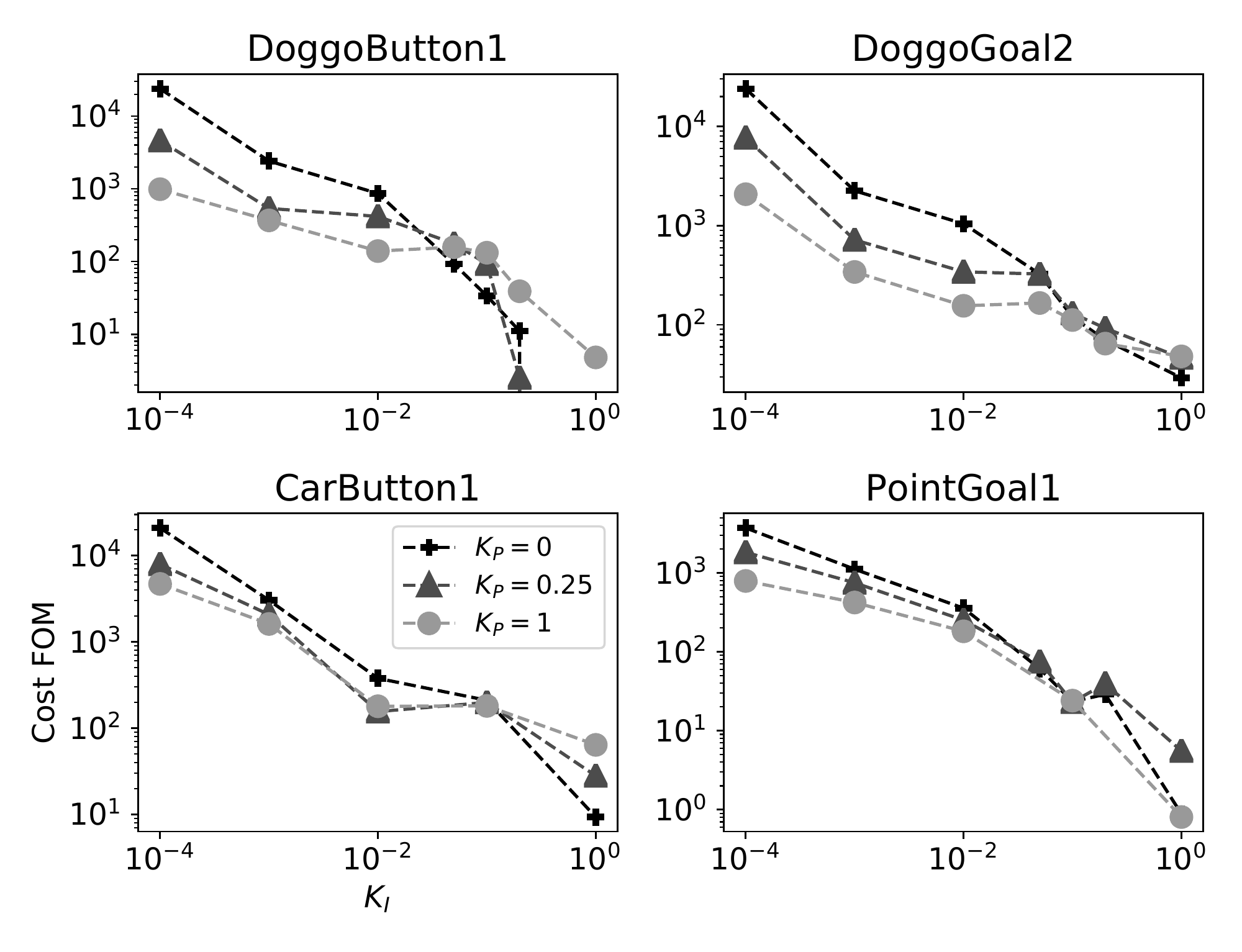}}
\vskip -0.2in
\caption{Learning run cost FOM versus penalty learning rate, $K_I$, from four environments spanning the robots in Safety Gym.  Each point is an average over four runs.  In all cases, PI-control improves performance (lower is better) over a wide and useful range of $K_I$, easing selection of that hyperparameter.}
\label{fig:fom_vs_alpha}
\end{center}
\vskip -0.3in
\end{figure}

We performed similar experiments on several Safety Gym environments in addition to \textsc{DoggoGoal2}: \textsc{PointGoal1}, the simplest domain with a point-like robot, \textsc{CarButton1}, for slightly more challenging locomotive control, and \textsc{DoggoButton1} for another challenging task (see appendix for learning curves like Figure \ref{fig:doggogoal_cost_return}).  Figure \ref{fig:fom_vs_alpha} plots the cost figure of merit over the same range of values for $K_I$, and for two strengths of added proportional control, for these environments.  PI-control clearly improved the cost FOM (lower is better) for $K_I<10^{-1}$, above which the fast integral control dominated.  Hence robustness to the value for $K_I$ was significantly improved in all the learning tasks studied.

\subsubsection{Control Efficiency}
\begin{figure}[ht]
\begin{center}
\centerline{
  \includegraphics[clip,width=0.8\columnwidth]{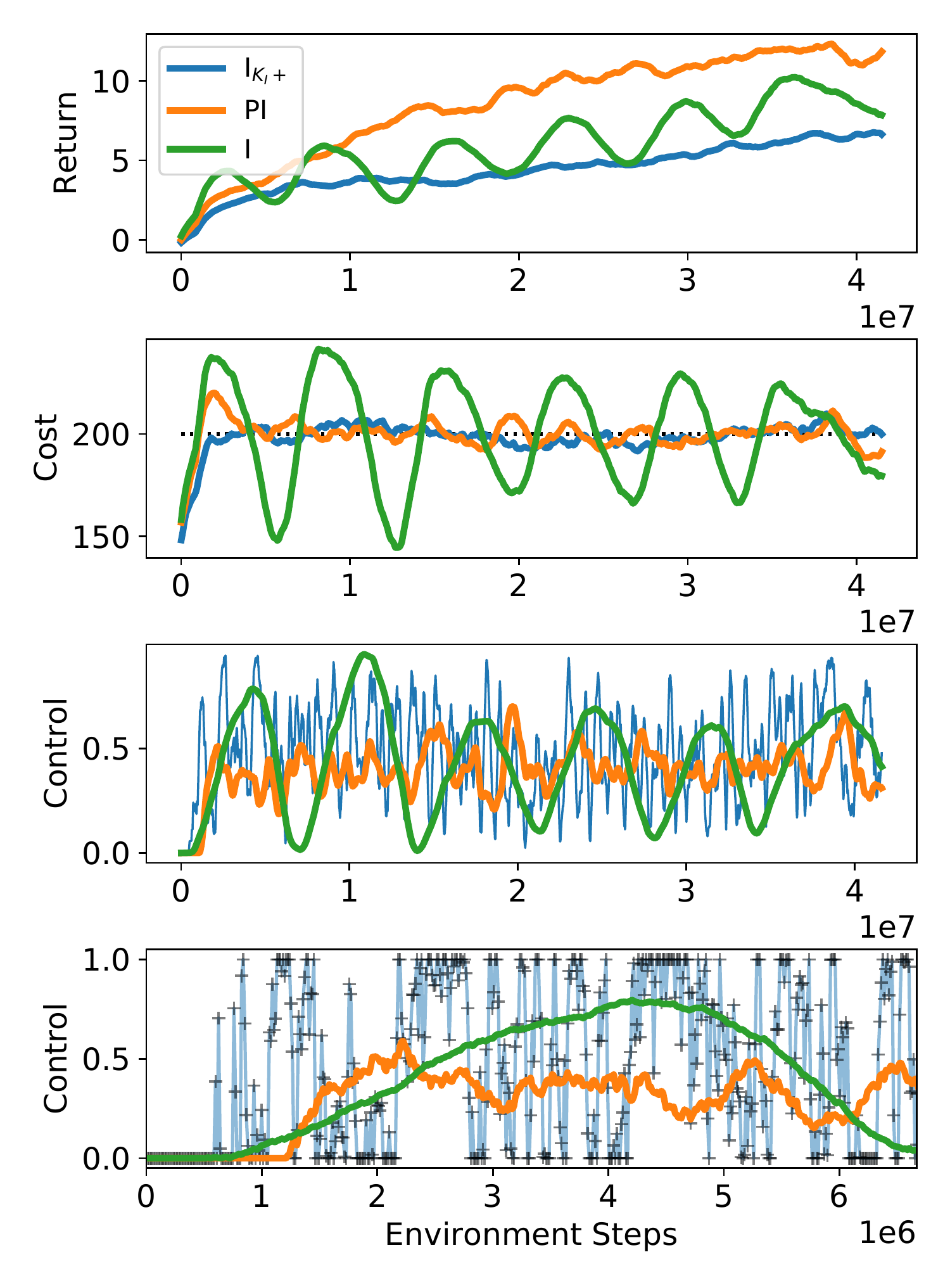}}
\caption{I- and PI-control with moderate $K_I=10^{-3}$ and I-control with fast $K_I=10^{-1}$ ($I_{K_I+}$).  \textit{Top} Returns diminished for fast-$K_I$, but high for PI. \textit{Second} Cost oscillations mostly damped by PI, removed by fast-$K_I$. \textit{Third} Control (smoothed) varies more rapidly under fast-$K_I$, is relatively steady for PI. \textit{Bottom} Control over first 500 RL iterations;  fast-$K_I$ slams the control to the extremes, causing the diminished returns.  Environment: \textsc{DoggoButton1}, cost limit 200.}
\label{fig:fast_a}
\end{center}
\vskip -0.2in
\end{figure}

We further investigated why increasing the penalty learning rate, $K_I$, eventually reduces reward performance, as was seen in the robustness study.  Figure \ref{fig:fast_a} shows learning curves for three settings: I- and PI-control with the same, moderate $K_I=10^{-3}$, and I-control with high $K_I=10^{-1}$.  The high-$K_I$ setting achieved responsive cost performance but lower long-term returns, which appears to result from  wildly fluctuating control.  In contrast, PI-control held relatively steady, despite the noise, allowing the agent to do reward-learning at every iteration.  The bottom panel displays individual control iterates, here displayed as $u=\lambda/(1+\lambda)$, over the first 7M environment steps, while the others show smoothed curves over the entire learning run, over 40M steps.

\subsubsection{Predictive Control by Derivatives}
\begin{figure}[ht]
\begin{center}
\centerline{
  \includegraphics[clip,width=1.0\columnwidth]{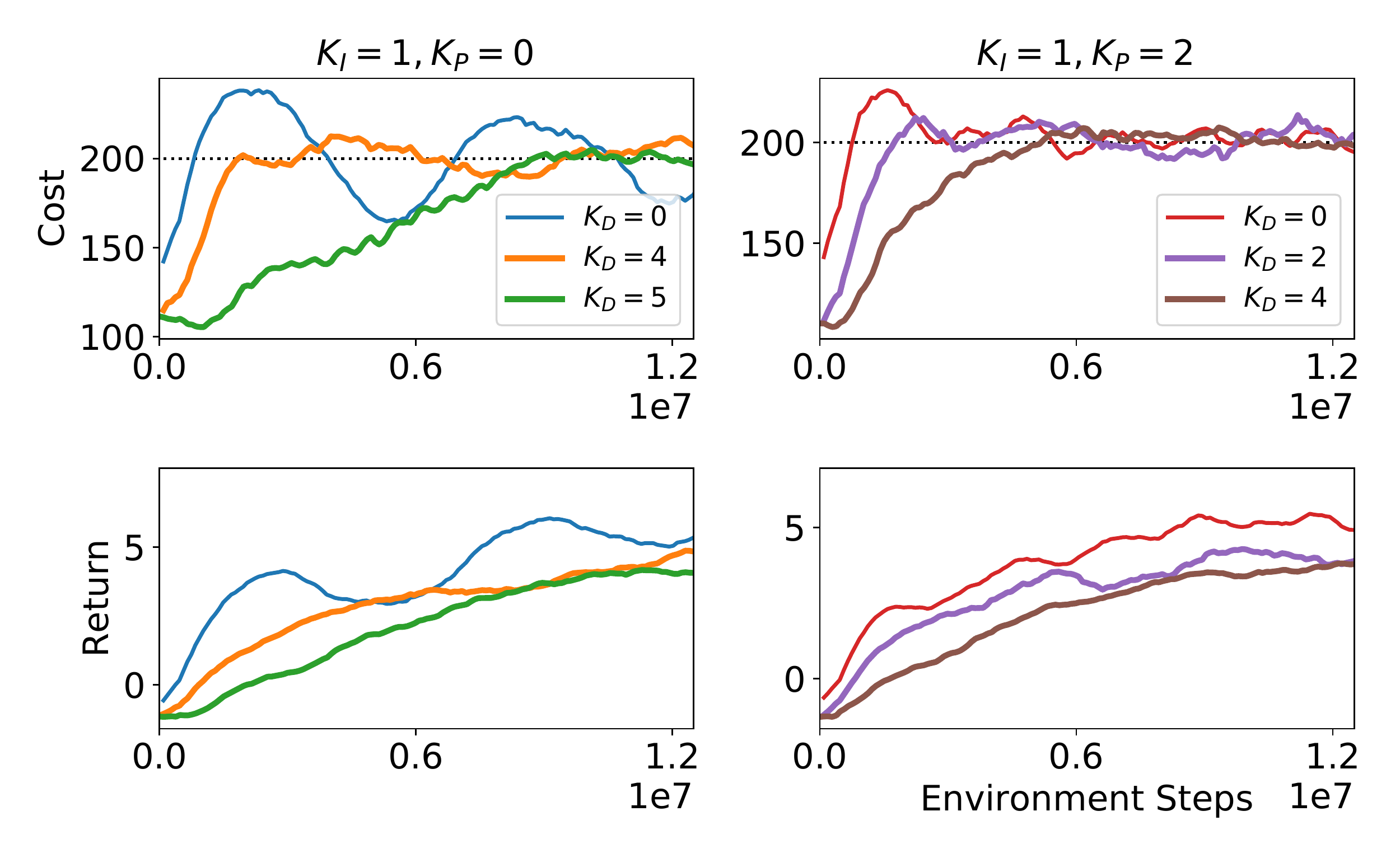}}
\vskip -0.1in
\caption{Derivative control can prevent cost overshoot and slow the rate of cost increase within feasible regions, which the Lagrangian method cannot do.  Environment: \textsc{DoggoButton1}, cost limit 200.}
\label{fig:derivative}
\end{center}
\vskip -0.2in
\end{figure}

Figure \ref{fig:derivative} demonstrates the predictive capabilities of derivative cost control in a noisy deep RL setting.  It removed cost overshoot from both the I- and PI-controlled baselines.  It was further able to slow the approach of the cost curve towards the limit, a desirable behavior for online learning systems requiring safety monitoring.  Curves for other environments are available in an appendix.

\section{Reward-Scale Invariance}
\label{sec:invariance}
In the preceding sections, we showed that PID control improves hyperparameter robustness in every constrained RL environment we tested.  Here we propose a complementary method to promote robustness both within and across environments.  Specifically, it addresses the sensitivity of learning dynamics to the relative numerical scale of reward and cost objectives.  

Consider two CMDPs that are identical except that in one the rewards are scaled by a constant factor, $\rho$.  The optimal policy parameters, $\theta^*$ remain unchanged, but clearly $\lambda^*$ must scale by $\rho$.  To attain the same learning dynamics, all controller settings, $\lambda_0, K_I, K_P$, and $K_D$ must therefore be scaled by $\rho$.  This situation might feature naturally within a collection of related learning environments.  Additionally, within the course of learning an individual CMDP, the balance between reward and cost magnitudes can change considerably, placing burden on the controller to track the necessary changes in the scale of $\lambda$.

One way to promote performance of a single choice of controller settings across these cases would be to maintain a fixed meaning for the value of $\lambda$ in terms of the relative influence of reward versus cost on the parameter update.  To this end, we introduce an adjustable scaling factor, $\beta_k$, in the policy gradient:
\begin{equation}
    \nabla_\theta\mathcal{L}=(1-u_k)\nabla_\theta J(\pi_{\theta_k}) - u_k\beta_k \nabla_\theta J_C(\pi_{\theta_k})
\end{equation}
A conspicuous choice for $\beta_k$ is the ratio of un-scaled policy gradients:
\begin{equation}
    \beta_{\nabla,k}=\frac{||\nabla_\theta J(\pi_{\theta_k})||}{||\nabla_\theta J_C(\pi_{\theta_k})||}
\end{equation}
since it balances the total gradient to have equal-magnitude contribution from reward- and cost-objectives at $\lambda=1$ and encourages $\lambda^*=1$. Furthermore, $\beta_\nabla$ is easily computed with existing algorithm components.  

\begin{figure}[ht]
\begin{center}
\centerline{
  \includegraphics[clip,width=1.\columnwidth]{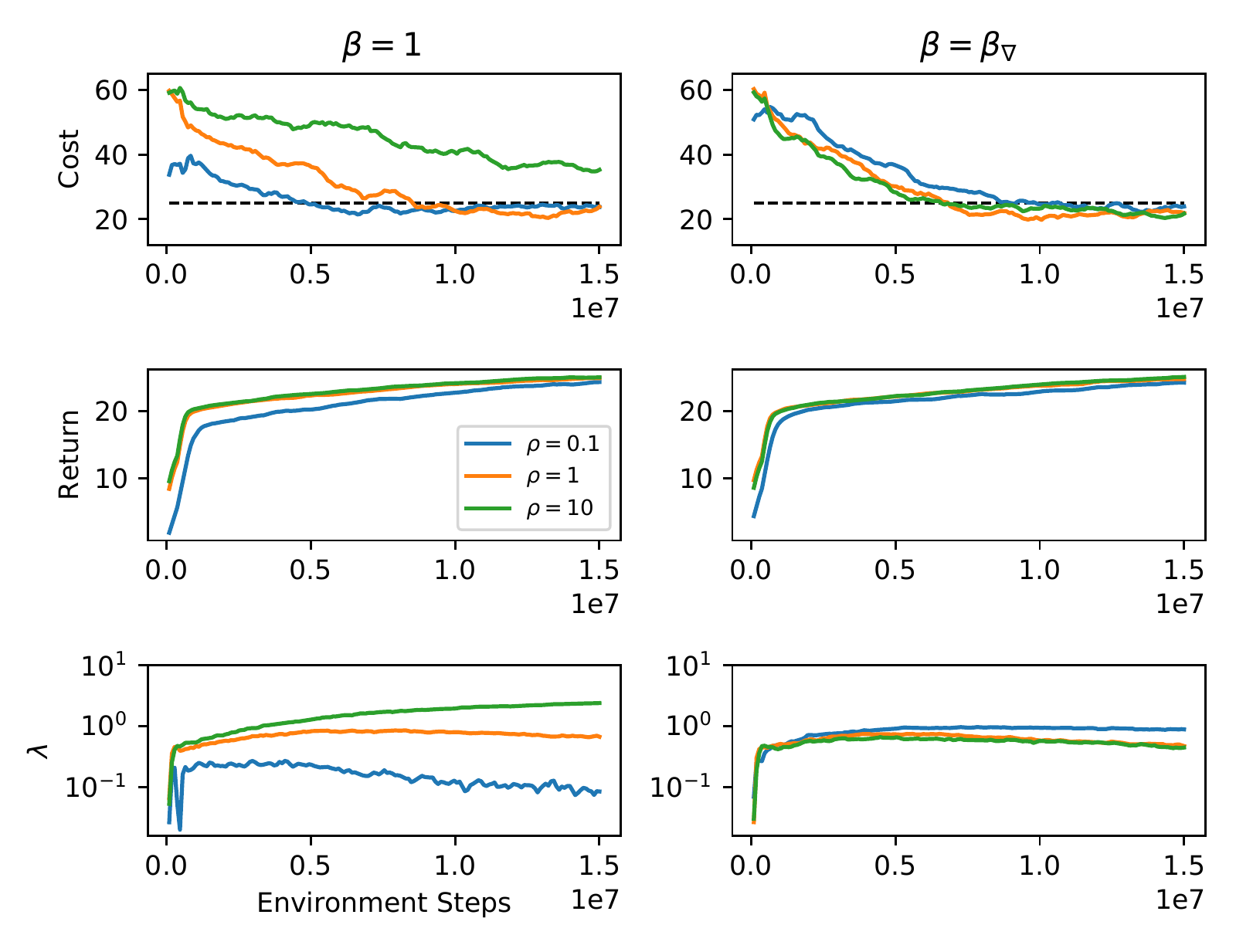}}
\vskip -0.1in
\caption{Costs, returns, and Lagrange multiplier with rewards scaled by $\rho\in\{0.1, 1, 10\}$; PI-control with $K_I=1e-3, K_P=0.1$.  \textit{Left column}: without objective-weighting, learning dynamics differ dramatically due to required scale of $\lambda$. \textit{Right column}: with objective-weighting, learning dynamics are nearly identical. Environment: \textsc{PointGoal1}, cost limit 25.}
\label{fig:reward_scale_PI}
\end{center}
\vskip -0.2in
\end{figure}

To test this method, we ran experiments on Safety Gym environments with their rewards scaled up or down by a factor of 10.  Figure \ref{fig:reward_scale_PI} shows a representative cross-section of results from the \textsc{PointGoal1} environment using PI-control.  The different curves within each plot correspond to different reward scaling.  Without objective-scaling (i.e. $\beta=1$), the dynamics under $\rho=10$ are as if controller parameters were instead divided by 10, and likewise for $\rho=0.1$.  Note the near-logarithmic spacing of $\lambda$ ($\lambda_{\rho=10}$ has not converged to its full value).  Using $\beta_\nabla$, on the other hand, the learning dynamics are nearly identical across two orders of magnitude of reward scale. $\lambda_0=1$ becomes an obvious choice for initialization, a point where previous theory provides little guidance \cite{chow2019} (although here we left $\lambda_0=0$).  Experiments in other environments and controller settings yielded similar results and are included in supplementary materials. Other methods, such running normalization of rewards and costs, could achieve similar effects and are worth investigating, but our simple technique is surprisingly effective and is not specific to RL.

\section{Conclusion}
\label{sec:conclusion}
Starting from a novel development in classic Lagrangian methods, we introduced a new set of constrained RL solutions which are straightforward to understand and implement, and we have shown them to be effective when paired with deep learning.

Several opportunities for further work lay ahead.  Analysis of the modified Lagrangian method and constrained RL as a dynamical system may relax theoretical requirements for a slowly-changing multiplier.  The mature field of control theory (and practice) provides tools for tuning controller parameters.  Lastly, the control-affine form may assist in both analysis (see \citet{how_much_uncertainty} and \citet{galbraith} for controllability properties for uncertain nonlinear dynamics) and by opening to further control techniques such as feedback linearization.  

Our contributions improve perhaps the most commonly used constrained RL algorithm, which is a workhorse baseline.  We have addressed its primary shortcoming while preserving its simplicity and even making it easier to use---a compelling combination to assist in a wide range of applications.

\section*{Acknowledgements}
Adam Stooke gratefully acknowledges previous support from the Fannie and John Hertz Foundation and the NVIDIA Corporation.  We thank Carlos Florensa and the anonymous reviewers for many helpful suggestions which improved the manuscript.

\bibliography{references}
\bibliographystyle{icml2020}

\end{document}


\onecolumn
\appendix

\icmltitle{Supplementary Materials for: Responsive Safety in Reinforcement Learning by PID Lagrangian Methods}

\section{Derivations}
Here we derive the effects of the proportional and derivative control terms on the Differential Multiplier Method.

\subsection{Traditional Lagrangian method}
Start from continuous Lagrangian dynamics for objective $f$ and constraint $g$:
\begin{align}
    \mathcal{L} &= f(\mathbf{x}) + \lambda g(\mathbf{x}) \\
    \dot{x}_i &= - \frac{\partial f}{\partial x_i} - \lambda \frac{\partial g}{\partial x_i} \\
    \dot{\lambda} &= \alpha g(\mathbf{x})
\end{align}
with $\alpha$ a scalar.  Taking the time-derivative of $\dot{x}_i$:
\begin{align}
    \ddot{x}_i &= -\frac{d}{dt}\left(\frac{\partial f}{\partial x_i}\right) - \dot{\lambda}\frac{\partial g}{\partial x_i} - \lambda \frac{d}{dt}\left(\frac{\partial g}{\partial x_i}\right) \\
    \frac{d}{dt}\left(\frac{\partial f}{\partial x_i}\right) &= \sum_j \frac{\partial}{\partial x_j} \left(\frac{\partial f}{\partial x_i}\right) \frac{dx_j}{dt} = \sum_j \frac{\partial^2 f}{\partial x_i \partial x_j}\dot{x}_j \\
    \frac{d}{dt}\left(\frac{\partial g}{\partial x_i}\right) &= \sum_j \frac{\partial^2 g}{\partial x_i \partial x_j}\dot{x}_j
\end{align}
and substituting for $\dot{\lambda}$:
\begin{align}
    \ddot{x}_i = -\sum_j \left(\frac{\partial^2 f}{\partial x_i \partial x_j} + \lambda \frac{\partial^2 g}{\partial x_i \partial x_j}\right)\dot{x}_j - \alpha g(\mathbf{x}) \frac{\partial g}{\partial x_i}
\end{align}
Defining the damping matrix, $A$:
\begin{align}
    A_{ij} &= \frac{\partial^2 f}{\partial x_i \partial x_j} + \lambda \frac{\partial^2 g}{\partial x_i \partial x_j}
\end{align}
yields the oscillator equation for the $i$-th component:
\begin{align}
    \ddot{x}_i + \sum_j A_{ij}\dot{x}_j + \alpha g(\mathbf{x}) \frac{\partial g}{\partial x_i} = 0
\end{align}
Written in vector form:
\begin{align}
    &\ddot{\mathbf{x}} + A\dot{\mathbf{x}} + \alpha g(\mathbf{x}) \nabla g = 0 \\
    &A = \nabla^2f + \lambda \nabla^2 g
\end{align}
with $\nabla^2 f$ denoting the Hessian of $f$.
\subsection{Proportional-Integral Method}
We amend the differential equation for $\lambda$ as:
\begin{align}
    \dot{\lambda} = \alpha g(\mathbf{x}) + \beta \dot{g}(\mathbf{x}) = \alpha g(\mathbf{x}) + \beta \sum_j \frac{\partial g}{\partial x_j}\dot{x}_j
\end{align}
Substitution into $\ddot{x}$ yields:
\begin{align}
    \ddot{x}_i &= -\sum_j \left(\frac{\partial^2 f}{\partial x_i \partial x_j} + \lambda \frac{\partial^2 g}{\partial x_i \partial x_j}\right)\dot{x}_j - \alpha g(\mathbf{x}) \frac{\partial g}{\partial x_i} - \beta \left(\sum_j\frac{\partial g}{\partial x_j}\dot{x}_j\right)\frac{\partial g}{\partial x_i} \\
    &= -\sum_j \left(\frac{\partial^2 f}{\partial x_i \partial x_j} + \lambda \frac{\partial^2 g}{\partial x_i \partial x_j} + \beta\frac{\partial g}{\partial x_i}\frac{\partial g}{\partial x_j}\right)\dot{x}_j - \alpha g(\mathbf{x}) \frac{\partial g}{\partial x_i}
\end{align}
giving the same dynamics as the traditional Lagrangian method except with damping matrix modified by addition of a positive-semi-definite term:
\begin{align}
        &\ddot{\mathbf{x}} + \left(A + \beta \nabla g \nabla^\top g \right)\dot{\mathbf{x}} + \alpha g(\mathbf{x}) \nabla g = 0
\end{align}

\subsection{Integral-Derivative Method}

We use the following differential equation for $\lambda$:
\begin{align}
    \dot{\lambda} = \alpha g(\mathbf{x}) + \gamma \ddot{g}(\mathbf{x})    
\end{align}
Expanding the second-derivative:
\begin{align}
    \ddot{g}(x) &= \frac{d}{dt}\left(\frac{d}{dt}g(\mathbf{x})\right) \\
    &= \frac{d}{dt}\left(\sum_j \frac{\partial g}{x_j}\dot{x_j}\right) \\
    &= \sum_j \left(\frac{d}{dt}\left(\frac{\partial g}{x_j}\right)\dot{x}_j + \frac{\partial g}{x_j}\frac{d}{dt}(\dot{x}_j)   \right) \\
    &= \sum_j \left( \dot{x}_j\sum_k\frac{\partial^2 g}{\partial x_k \partial x_j}\dot{x}_k + \frac{\partial g}{\partial x_j}\ddot{x}_j  \right)
\end{align}
Substituting the full expression for $\dot{\lambda}$ into $\ddot{x_i}$ yields:
\begin{align}
    \ddot{x}_i = -\sum_j \left( \frac{\partial^2 f}{\partial x_i \partial x_j} + \lambda \frac{\partial^2 g}{\partial x_i \partial x_j} \right)\dot{x_j} - \alpha g(\mathbf{x})\frac{\partial g}{\partial x_i} - \gamma \frac{\partial g}{\partial x_i} \left( \sum_j\left(\sum_k \dot{x}_j \frac{\partial^2 g}{\partial x_j \partial x_k} \dot{x}_k + \frac{\partial g}{\partial x_j}\ddot{x}_j\right) \right)
\end{align}
where the terms with coefficient $\gamma$ are due to the derivative update term in $\lambda$.  Now included are terms of second-order in the velocity and additional mixing terms on the acceleration.  It is instructive to consider the resulting vector equation in full:
\begin{align}
    0 &= \ddot{\mathbf{x}} + A\dot{\mathbf{x}} + \alpha g(\mathbf{x}) \nabla g + \gamma \nabla g \dot{\mathbf{x}}^\top \nabla^2 g \dot{\mathbf{x}} + \gamma \nabla g \nabla^\top g \ddot{\mathbf{x}} \\
    &= (I + \gamma \nabla g \nabla^\top g)\ddot{\mathbf{x}} + A\dot{\mathbf{x}} + \left(\alpha g(\mathbf{x}) + \gamma \dot{\mathbf{x}}^\top \nabla^2 g \dot{\mathbf{x}}  \right)  \nabla g 
\end{align}
The acceleration terms are coupled together by a positive definite matrix (identity plus the outer product of the vector $\nabla g$ with itself).  Therefore a form without coupling of acceleration can be restored by left-multiplying by the inverse of this matrix, letting $B=(I + \gamma \nabla g \nabla^\top g)$:
\begin{align}
    \ddot{\mathbf{x}} + B^{-1}A\dot{\mathbf{x}} + \left(\alpha g(\mathbf{x}) + \gamma \dot{\mathbf{x}}^\top \nabla^2 g \dot{\mathbf{x}}  \right)  B^{-1} \nabla g = 0
\end{align}
The effects of thenew terms from the derivative update rule are discussed in the main text.

\subsection{Proportional-Integral-Derivative Method}
Finally, the full PID-Lagrangian method has dynamics which combines the independent effects of the proportional and derivative methods:
\begin{align}
    \ddot{\mathbf{x}} + B^{-1}\left(A + \beta \nabla g \nabla^\top g \right)\dot{\mathbf{x}} + \left(\alpha g(\mathbf{x}) + \gamma \dot{\mathbf{x}}^\top \nabla^2 g \dot{\mathbf{x}}  \right)  B^{-1} \nabla g = 0
\end{align}

\section{Training Details}
All of our experiments began with randomly initialized agents.  The reward-value and cost-value estimators shared parameters with the policy.  We used Generalized Advantage Estimation for both reward and cost advantages.  The control input is updated once per iteration, which in our settings included multiple gradient updates on the policy.  Training batches typically included the end of several trajectories, allowing an estimate of the average episodic sum of costs at each iteration.

\begin{table}[h]
\caption{Experiment hyperparameters.}
\label{table:params}
\vskip 0.15in
\begin{center}
\begin{small}
\begin{sc}
\begin{tabular}{rl}
\toprule
Hyperparameter & Value \\
\midrule
learning rate & $1\times10^{-4}$ \\
NN hidden layer size & 512 \\
NN nonlinearity & $tanh$ \\
batch dimension, time & 128 \\
batch dimension, envs & 104 \\
PPO epochs & 8 \\
PPO minibatches & 1 \\
PPO ratio-clip & 0.1 \\
Discount, $\gamma$ & 0.99 \\
$\lambda_{GAE}$ & 0.97 \\
Cost scaling & $1/10$ \\
Exponential moving average, $K_P$ & 0.95 \\
Exponential moving average, $K_D$ & 0.9 \\
Difference iterates delay, $K_D$ & 15 \\
Observation normalization & True \\
$\hat{\beta}$ & 1 (unless specified) \\
Exponential moving average, $\hat{\beta}$ & 0.9 \\
\bottomrule
\end{tabular}
\end{sc}
\end{small}
\end{center}
\vskip -0.1in
\end{table}

\newpage
\section{Additional Learning Curves}
\subsection{Cost and Reward Curves for Additional Environments}
This section contains learning curves from the experiments used to make the figures showing cost figure-of-merit versus Lagrange multiplier learning rate, $K_I$, in the main text.  The main text includes curves from \textsc{DoggoGoal2}, the other environments are shown here.

\begin{figure}[ht]
\begin{center}
\centerline{
  \includegraphics[clip,width=0.7\columnwidth]{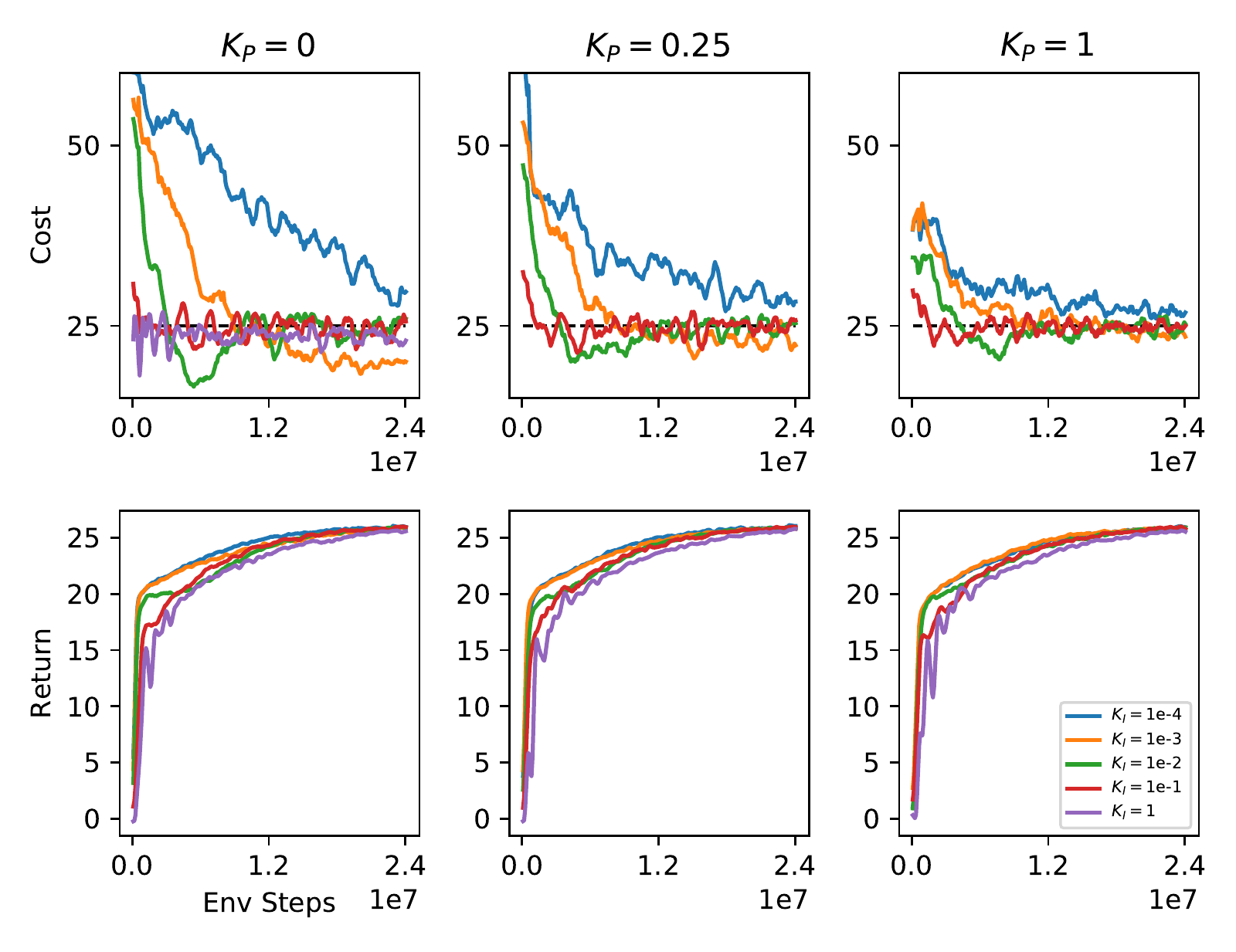}}
\caption{Costs and returns with varying Lagrange multiplier learning rate, $K_I$, and proportional control coefficient, $K_P$, in \textsc{PointGoal1}, cost-limit=25.}
\label{fig:cost_returns_pointgoal}
\end{center}
\end{figure}

\begin{figure}[ht]
\begin{center}
\centerline{
  \includegraphics[clip,width=0.7\columnwidth]{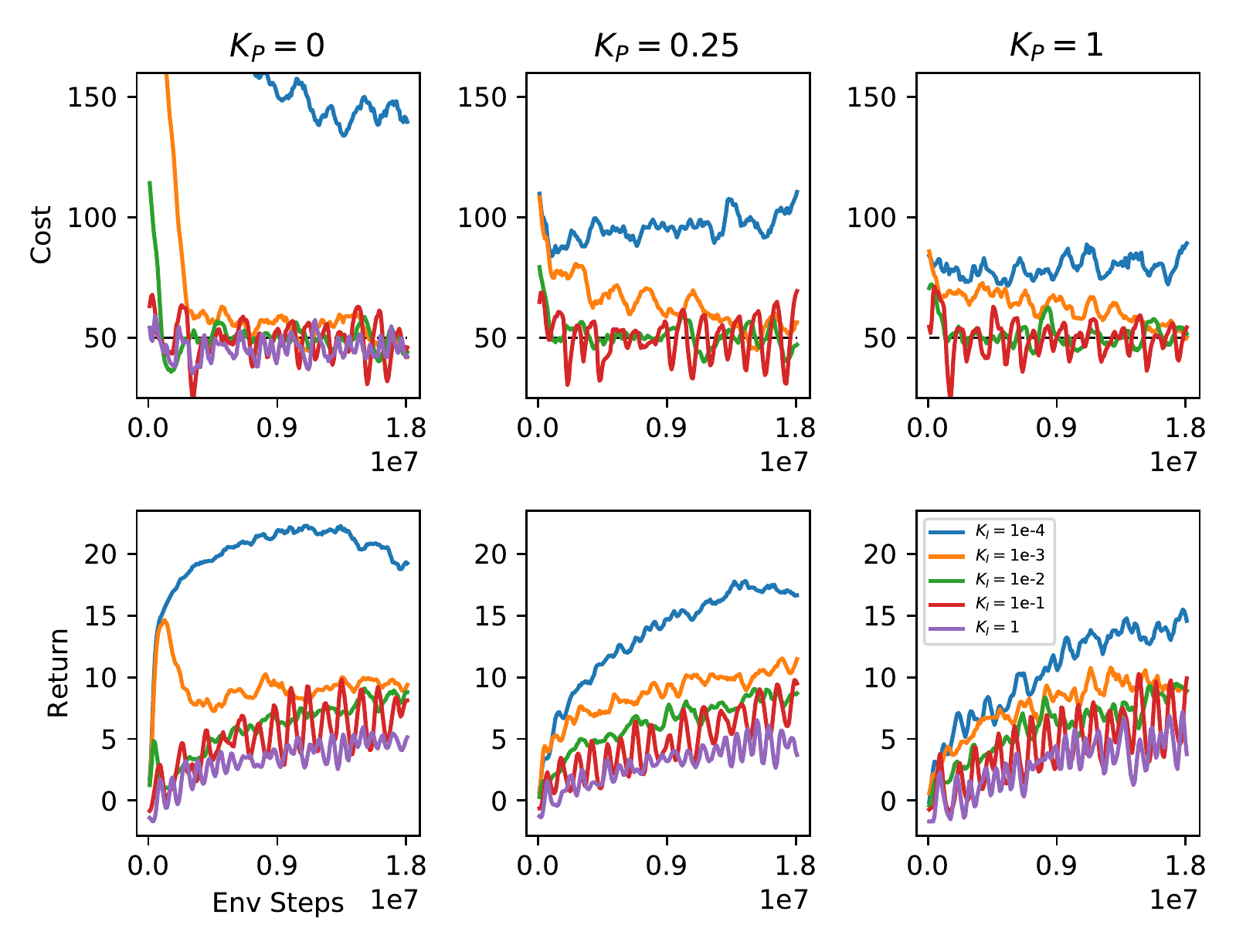}}
\caption{Costs and returns with varying Lagrange multiplier learning rate, $K_I$, and proportional control coefficient, $K_P$, in \textsc{CarButton1}, cost-limit=50.}
\label{fig:cost_returns_carbutton}
\end{center}
\end{figure}

\begin{figure}[ht]
\begin{center}
\centerline{
  \includegraphics[clip,width=0.7\columnwidth]{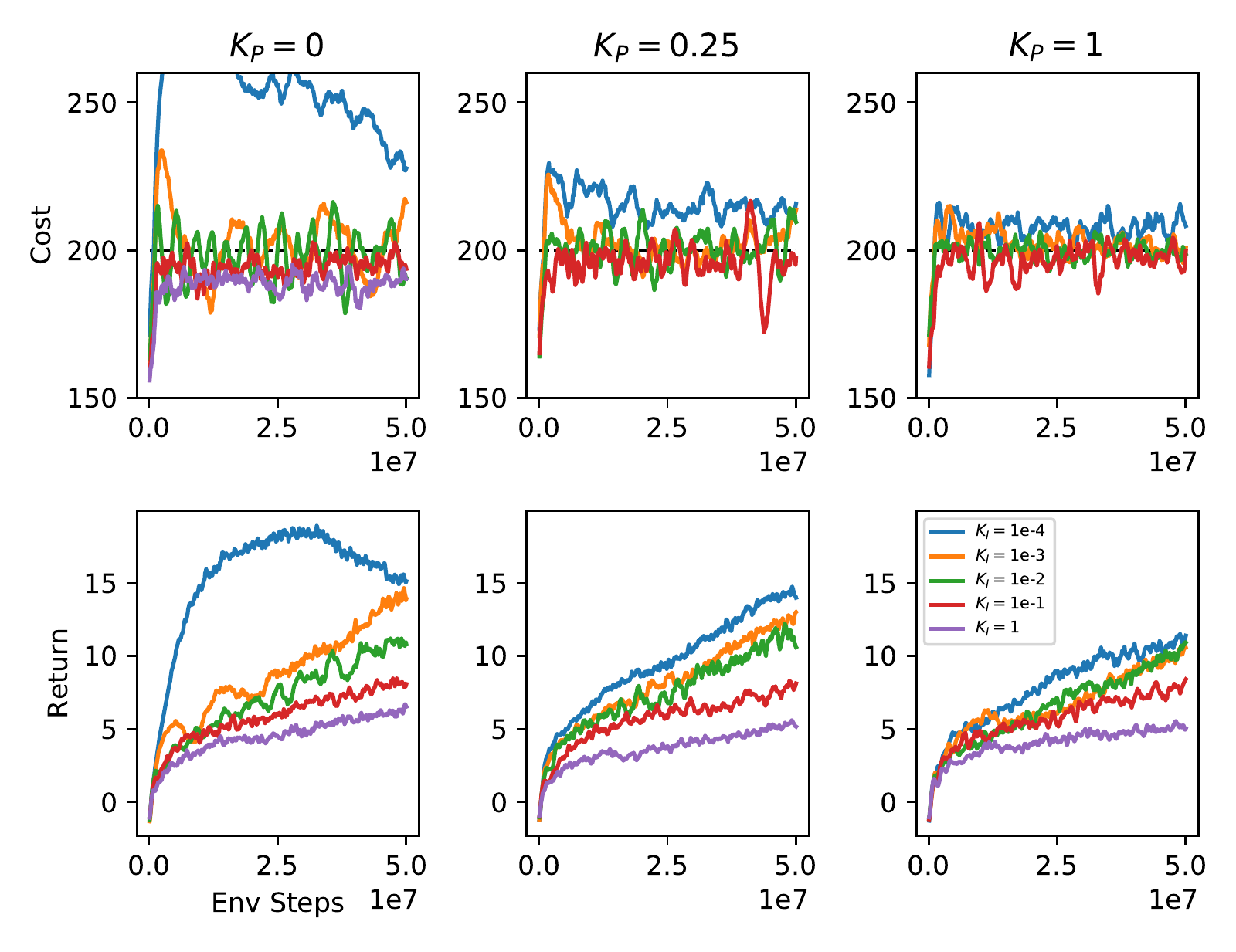}}
\caption{Costs and returns with varying Lagrange multiplier learning rate, $K_I$, and proportional control coefficient, $K_P$, in \textsc{DoggoButton1}, cost-limit=200.}
\label{fig:cost_returns_doggobutton}
\end{center}
\end{figure}

\clearpage
\subsection{Derivative-Control Learning Curves}

\begin{figure}[h]
\centering
\begin{minipage}{.45\textwidth}
  \centering
  \includegraphics[width=.7\linewidth]{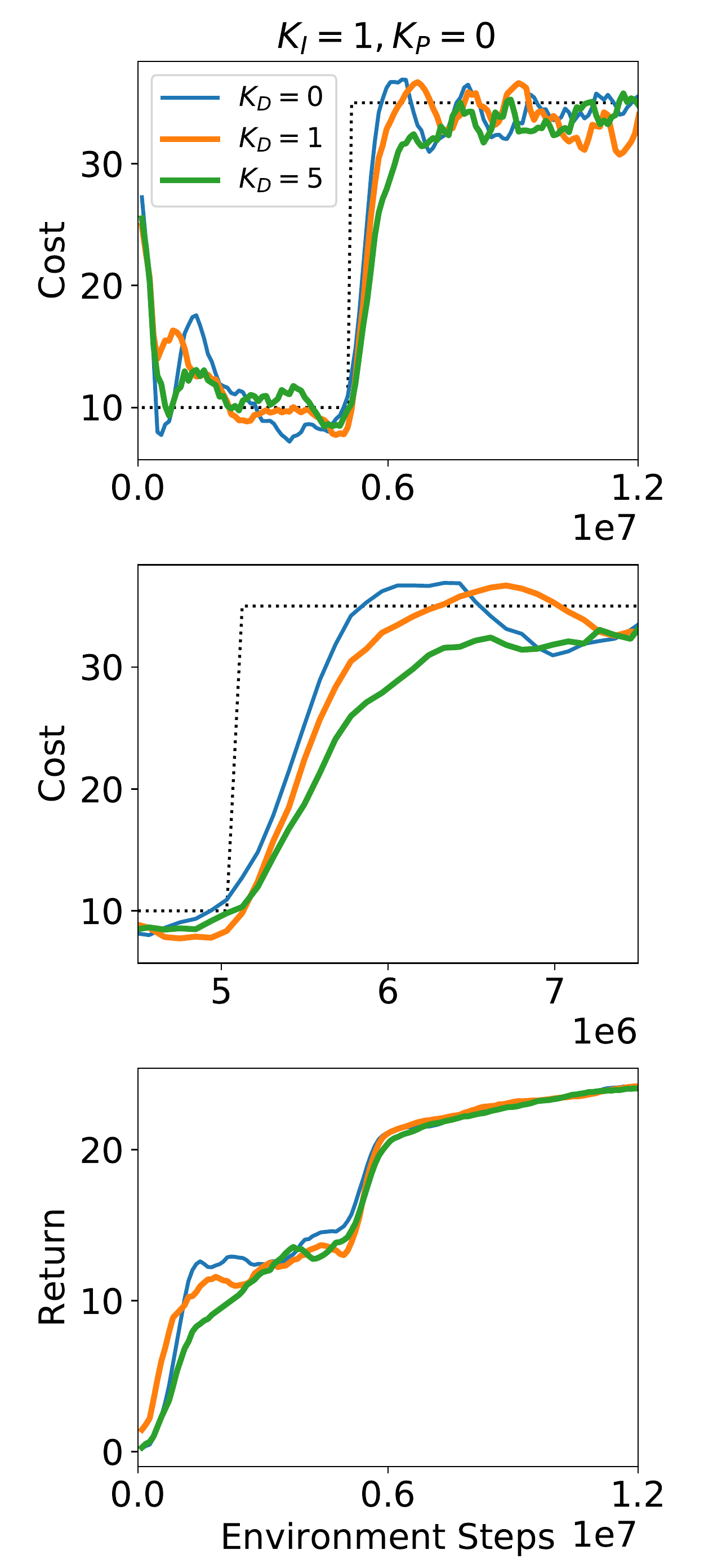}
  \captionof{figure}{Derivative control slows the increase in cost, causing it to rise more gradually, and reducing overshoot.  The cost limit was increased from 10 to 35 at 5M environment steps.  Environment: \textsc{PointGoal1}.}
  \label{fig:test1}
\end{minipage}%
\hspace{12pt}
\begin{minipage}{.45\textwidth}
  \centering
  \includegraphics[width=.7\linewidth]{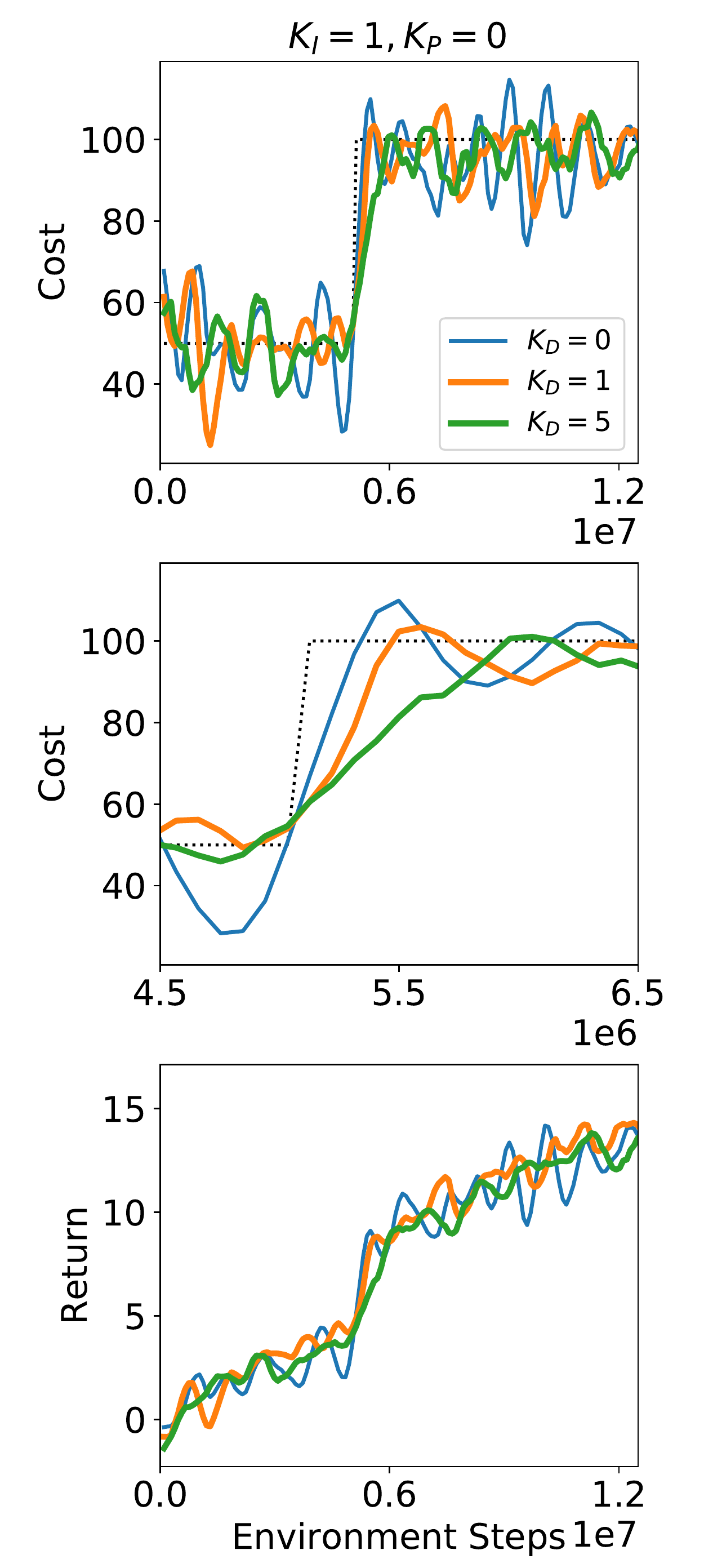}
  \captionof{figure}{Derivative control slows the increase in cost, causing it to rise more gradually, and reducing overshoot.  The cost limit was increased from 50 to 100 at 5M environment steps.  Environment: \textsc{CarButton1}.}
  \label{fig:test2}
\end{minipage}
\end{figure}



\clearpage
\subsection{Comparison to Unconstrained Algorithms}
This section shows learning curves for the same four environments referenced in the main text.  These figures demonstrate that the unconstrained algorithm (PPO) does not satisfy the cost constraints, and as a result it achieves higher rewards.

\begin{figure}[ht]
  \begin{subfigure}
    \centering
    \includegraphics[width=0.45\columnwidth]{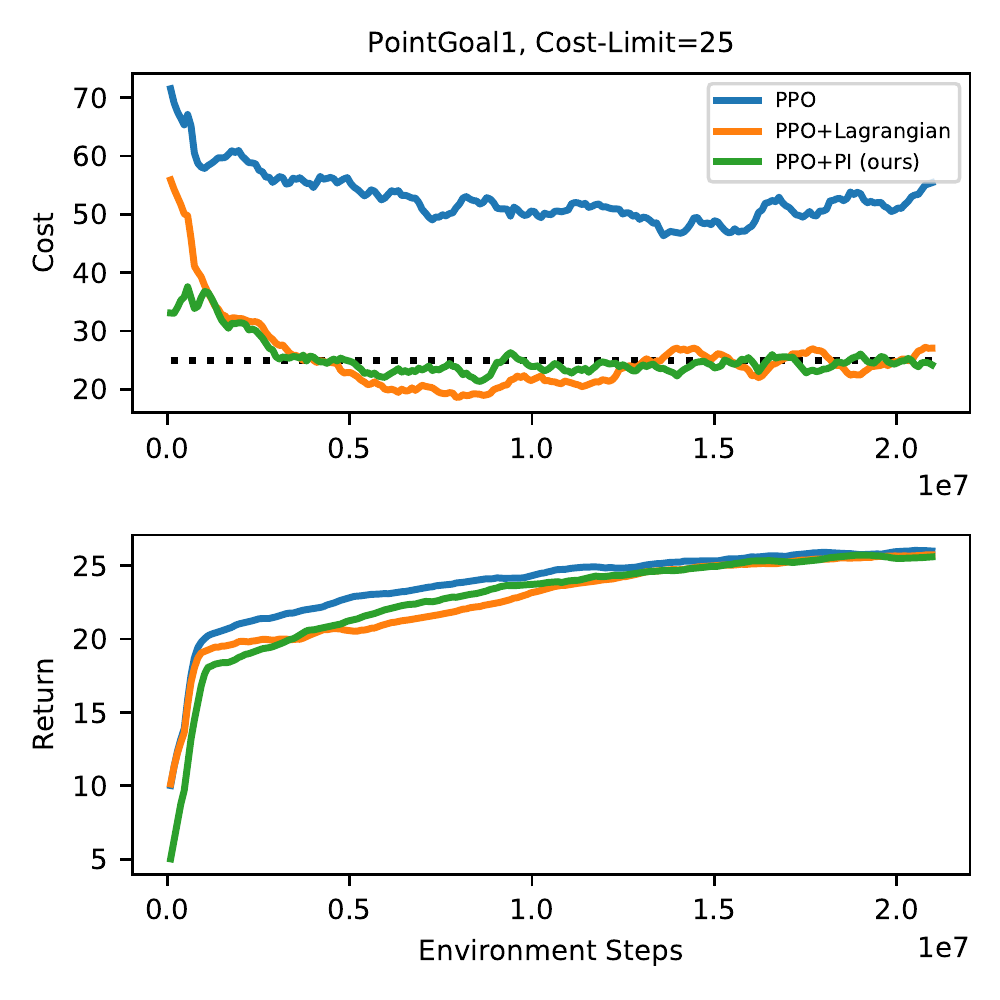}
  \end{subfigure}
  \hfill
  \begin{subfigure}
    \centering
    \includegraphics[width=0.45\columnwidth]{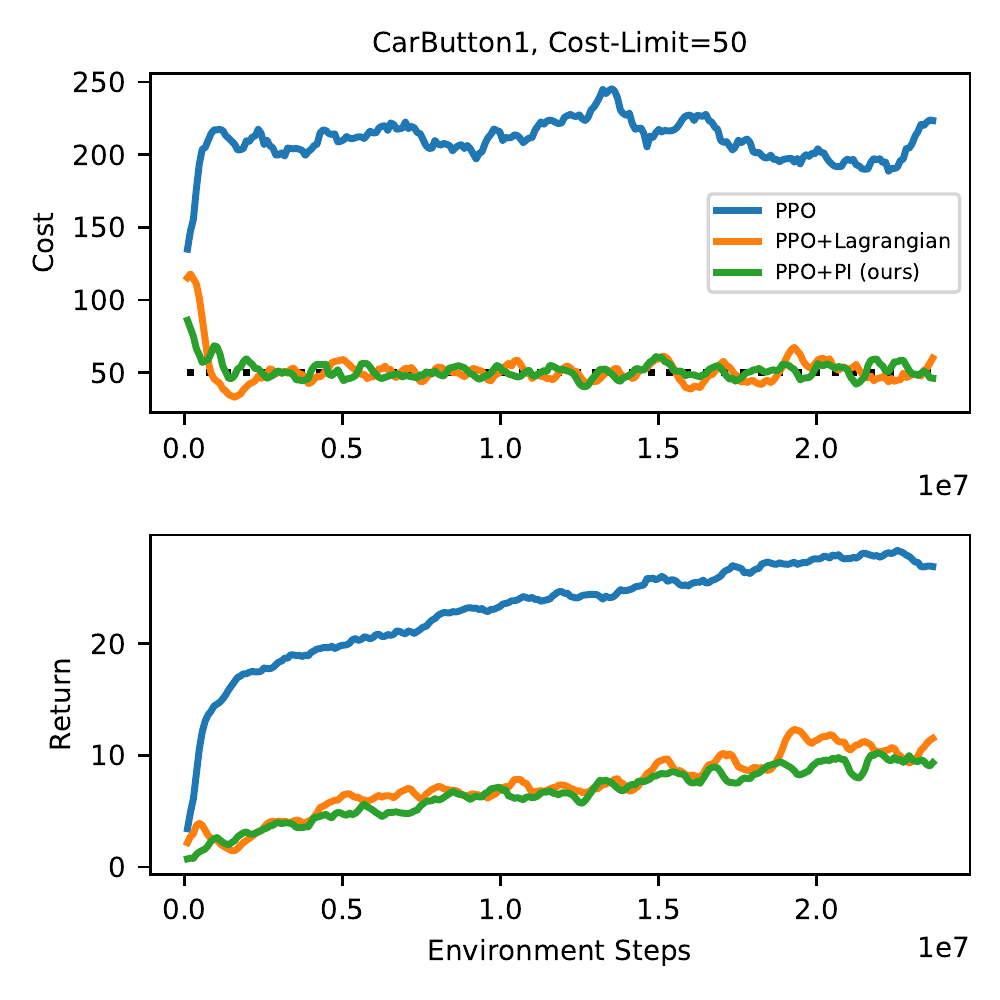}
  \end{subfigure}

  \medskip

  \begin{subfigure}
    \centering
    \includegraphics[width=0.45\columnwidth]{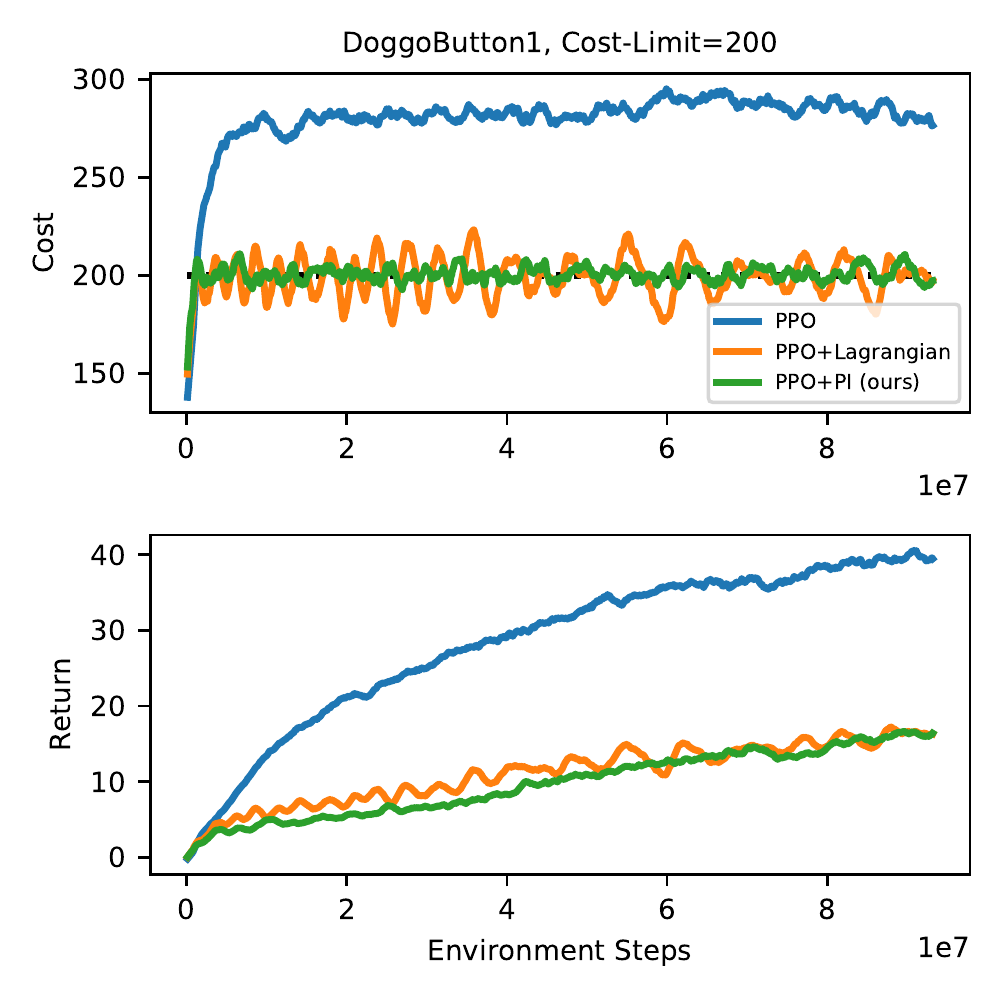}
  \end{subfigure}
  \hfill
  \begin{subfigure}
    \centering
    \includegraphics[width=0.45\columnwidth]{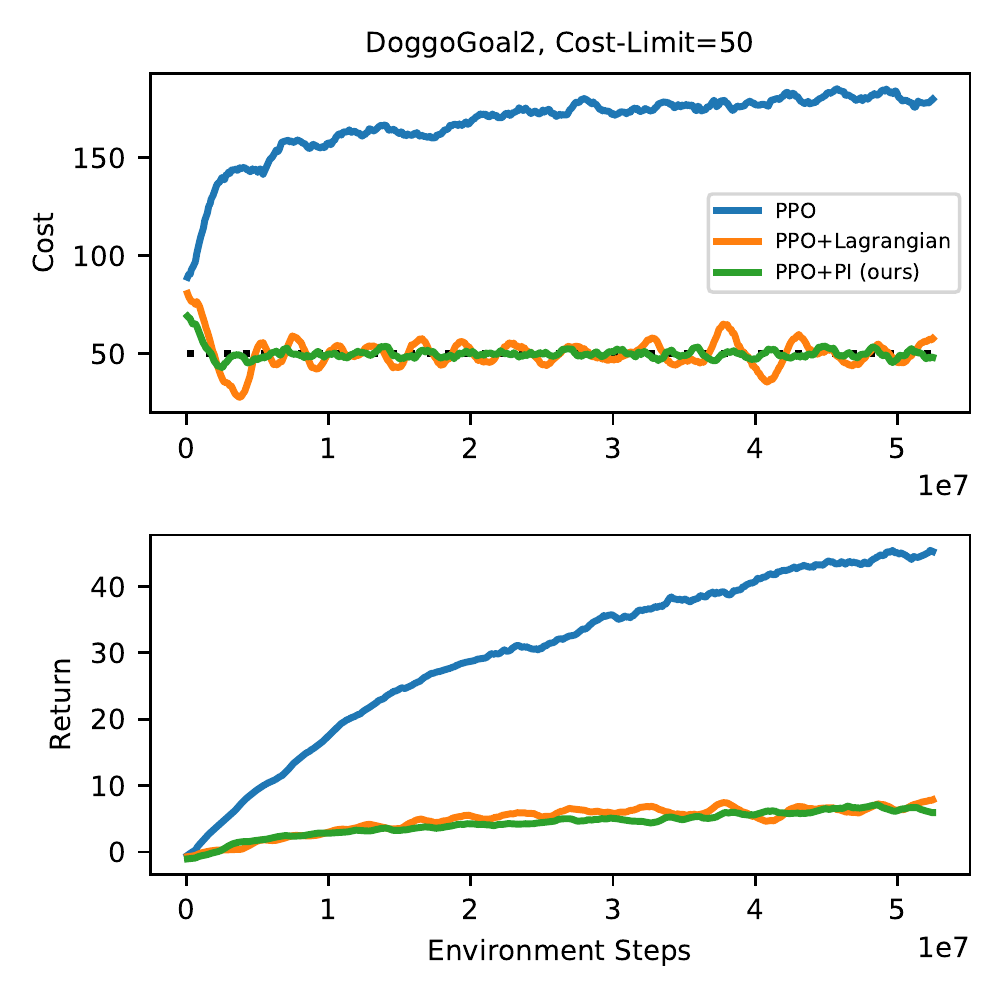}
  \end{subfigure}
\caption{Cost and reward curves for three variants of PPO: unconstrained , Lagrangian, and PI-Controlled.  The unconstrained algorithm wildly violates all cost-limits used in our experiments.  PPO+Lagrangian and PPO+PI use the same Lagrange multiplier learning rate, $K_I$.}
\label{fig:unconstrained}
\end{figure}

\clearpage
\section{Adaptive Objective-Balancing}
\textbf{Alternative, KL-Based Estimator}
The magnitudes of the gradients in $\theta$-space might not fully reflect the relative impacts of reward- and cost-learning on agent behavior.  Reinforcement learning offers an alternative grounding in policy-space, for example using:
\begin{equation}
    \beta_{KL}=\frac{D_{KL}(\pi_{\theta_k}||\pi_{\theta_k+1}^R)}{D_{KL}(\pi_{\theta_k}||\pi_{\theta_k+1}^C)}
\end{equation}
Here, $\pi_{\theta_{k+1}}^R$ and $\pi_{\theta_{k+1}}^C$ are hypothetical new policies found using only the reward-objective or only the (un-scaled) cost-objective, respectively.  We experimented with this estimator and found it to work, although not as well as the gradient-norm estimator in our cases.  Some results are included in the figures below.

\textbf{Figures}
We include figures for \textsc{PointGoal1} using I-control and PI-control, and the more challenging \textsc{DoggoGoal2} using I-control.  Observe in the un-balanced case ($\hat{\beta}=1$) that as the controller settings scale with the reward, the learning dynamics remained the same.  For example, see
$(K_I=0.1, \rho=10), (K_I=0.01, \rho=1),\text{ and } (K_I=0.001, \rho=0.1)$.  Using gradient-based objective balancing ($\beta_\nabla$), however, the learning dynamics were roughly the same across all reward scales, for given controller settings.  Alternatively, KL-based objective balancing ($\beta_{KL}$) was also effective, but did not produce dynamics as uniform as the gradient-based method.

\begin{figure}[ht]
\begin{center}
\underline{$K_I=0.1$}
\centerline{
  \includegraphics[clip,width=0.8\columnwidth]{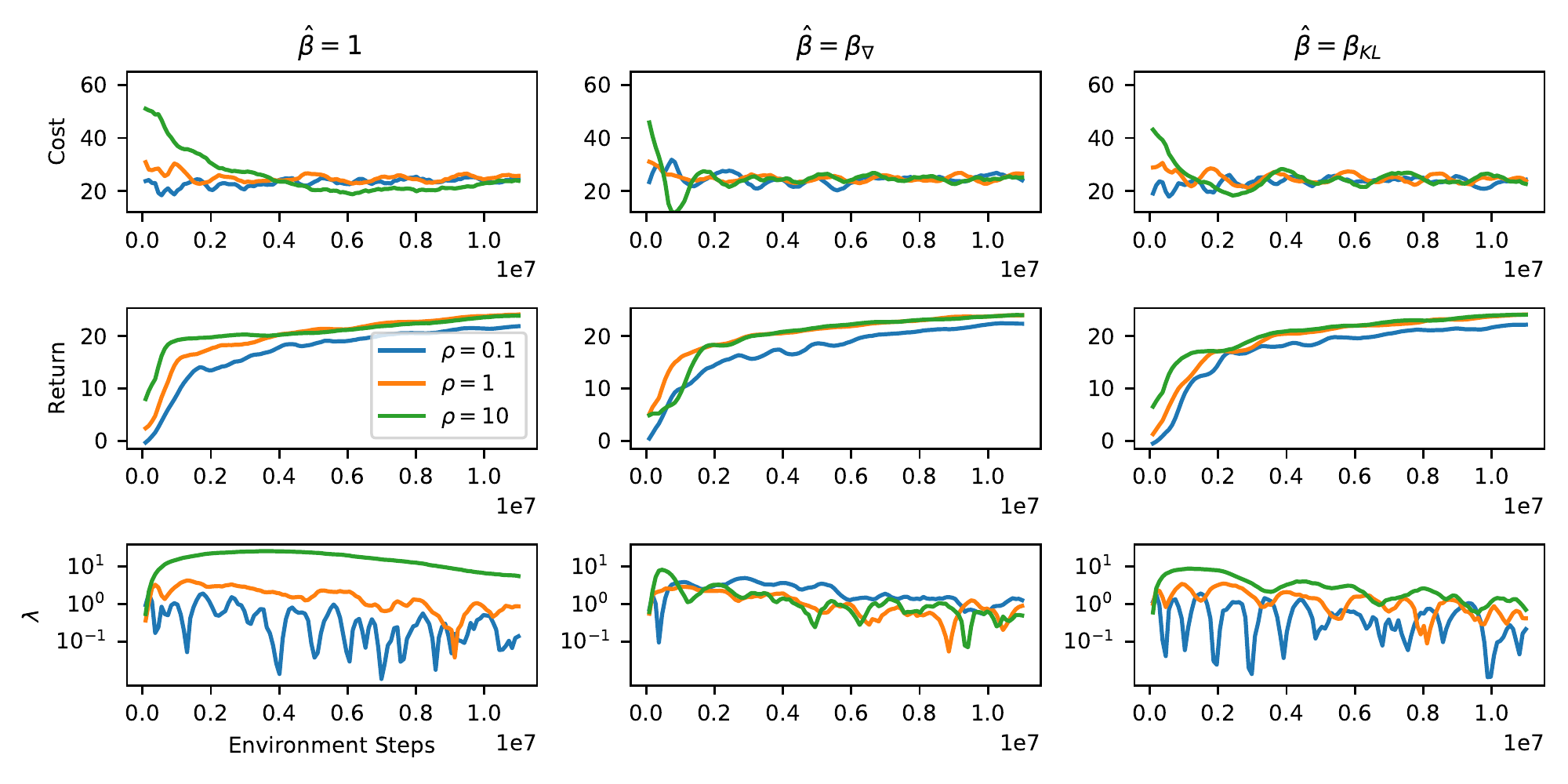}}
\underline{$K_I=0.01$}
\centerline{
  \includegraphics[clip,width=0.8\columnwidth]{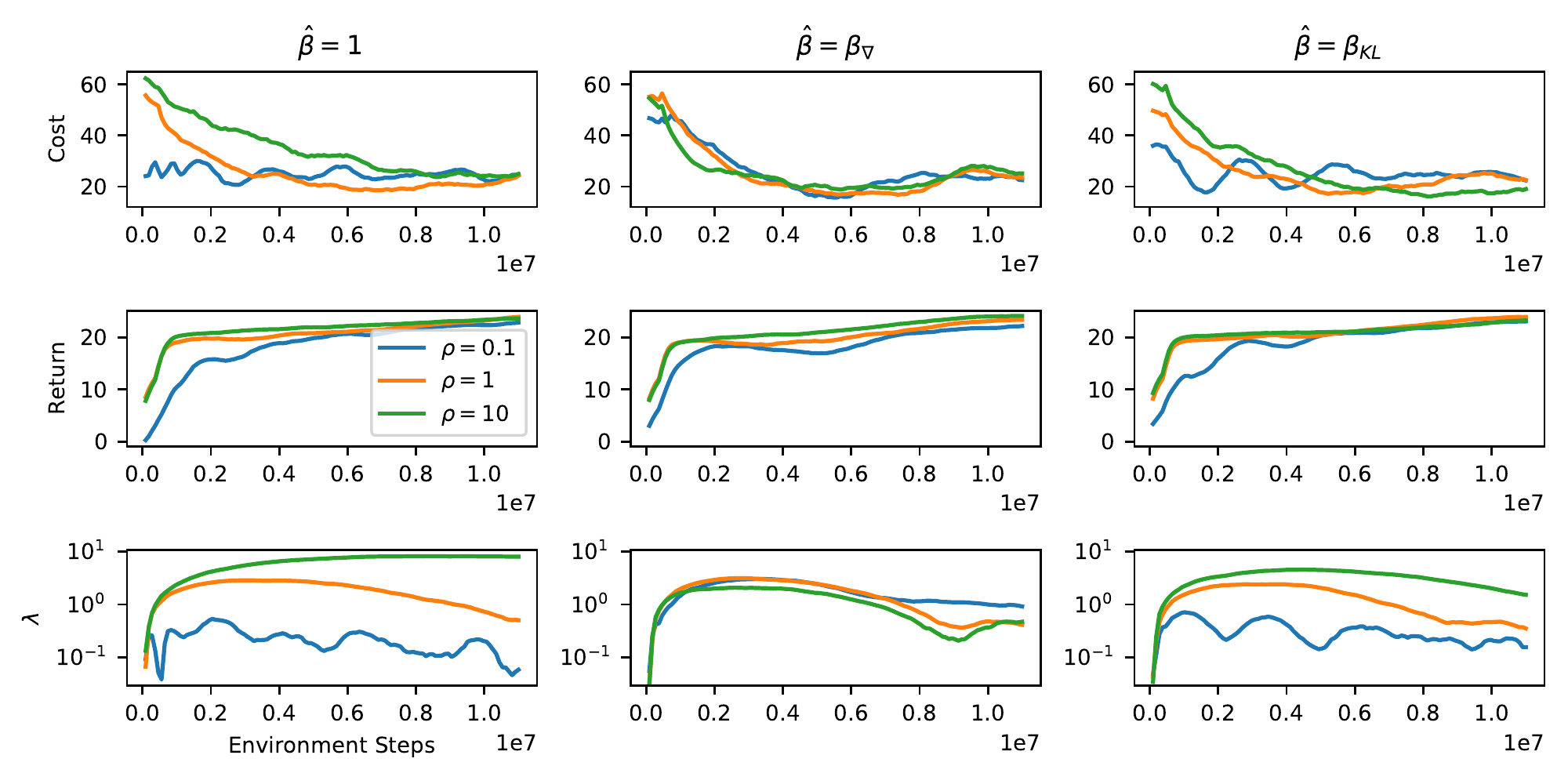}}
\underline{$K_I=0.001$}
\centerline{
  \includegraphics[clip,width=0.8\columnwidth]{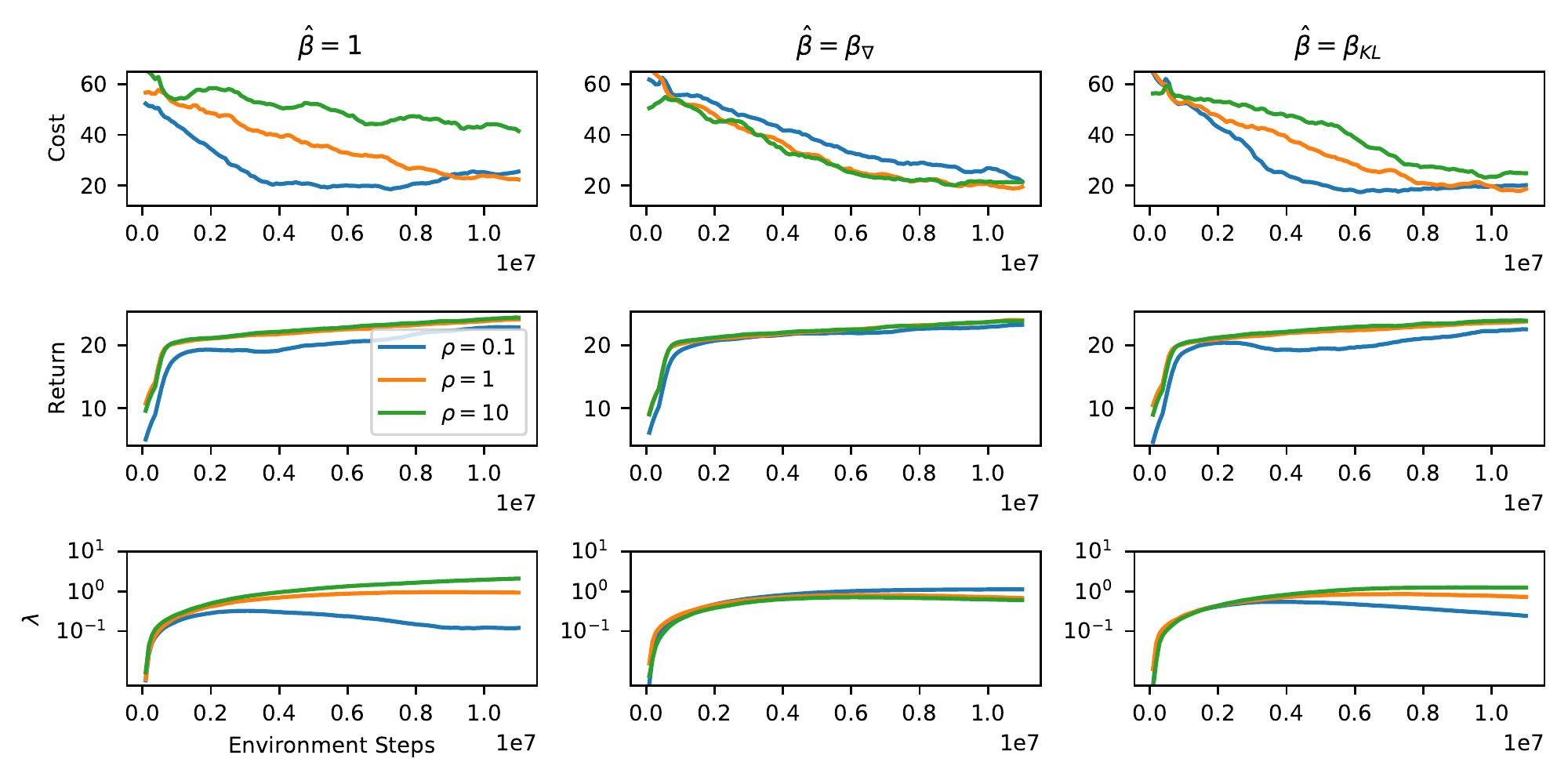}}
\caption{Reward scaling, I-control, \textsc{PointGoal1}, cost-limit=25.}
\label{fig:rewardscale_pointgoal_alpha}
\end{center}
\end{figure}

\begin{figure}[ht]
\begin{center}
\underline{$K_P=0.1$}
\centerline{
  \includegraphics[clip,width=0.8\columnwidth]{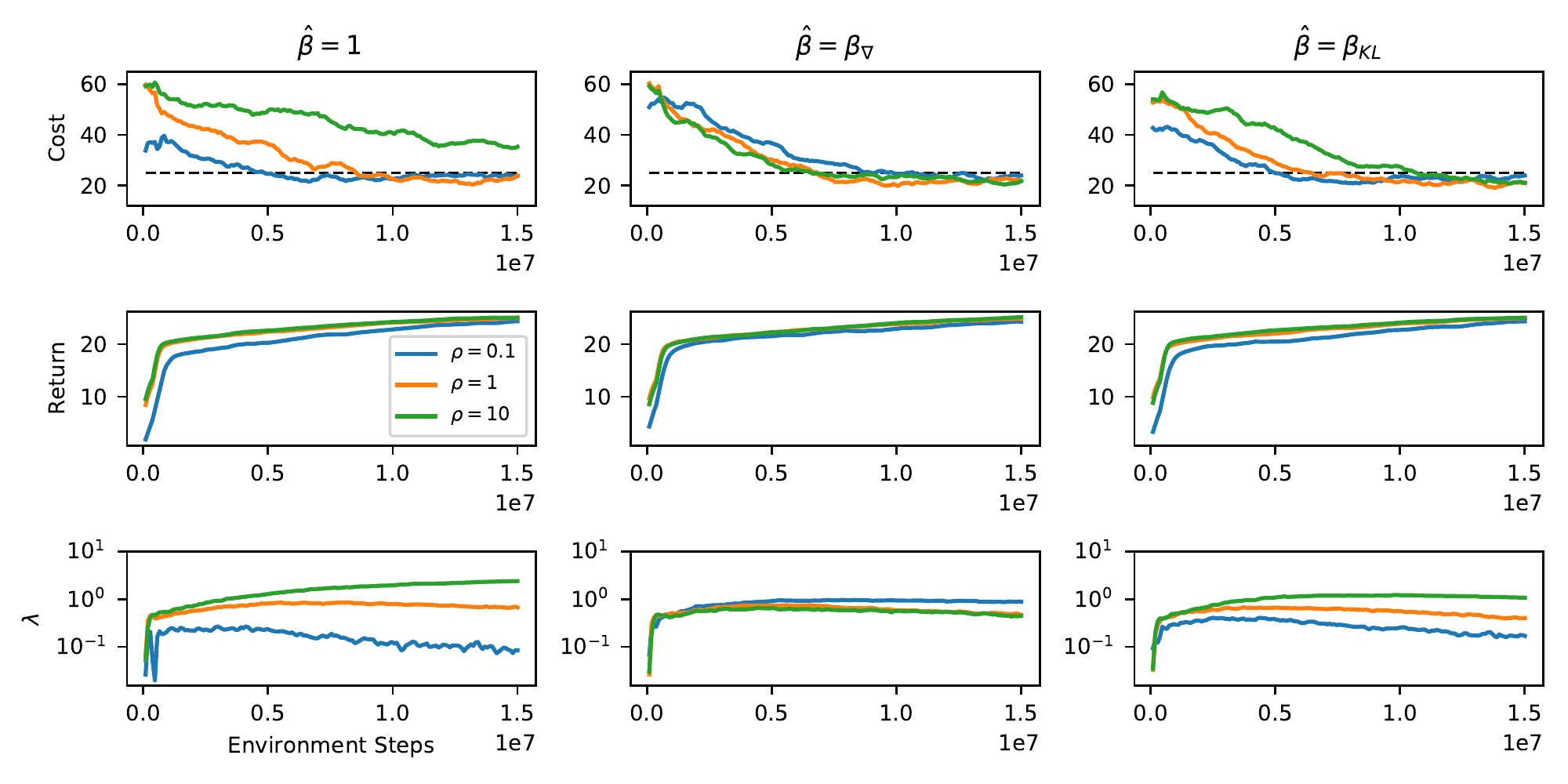}}
\underline{$K_P=1$}
\centerline{
  \includegraphics[clip,width=0.8\columnwidth]{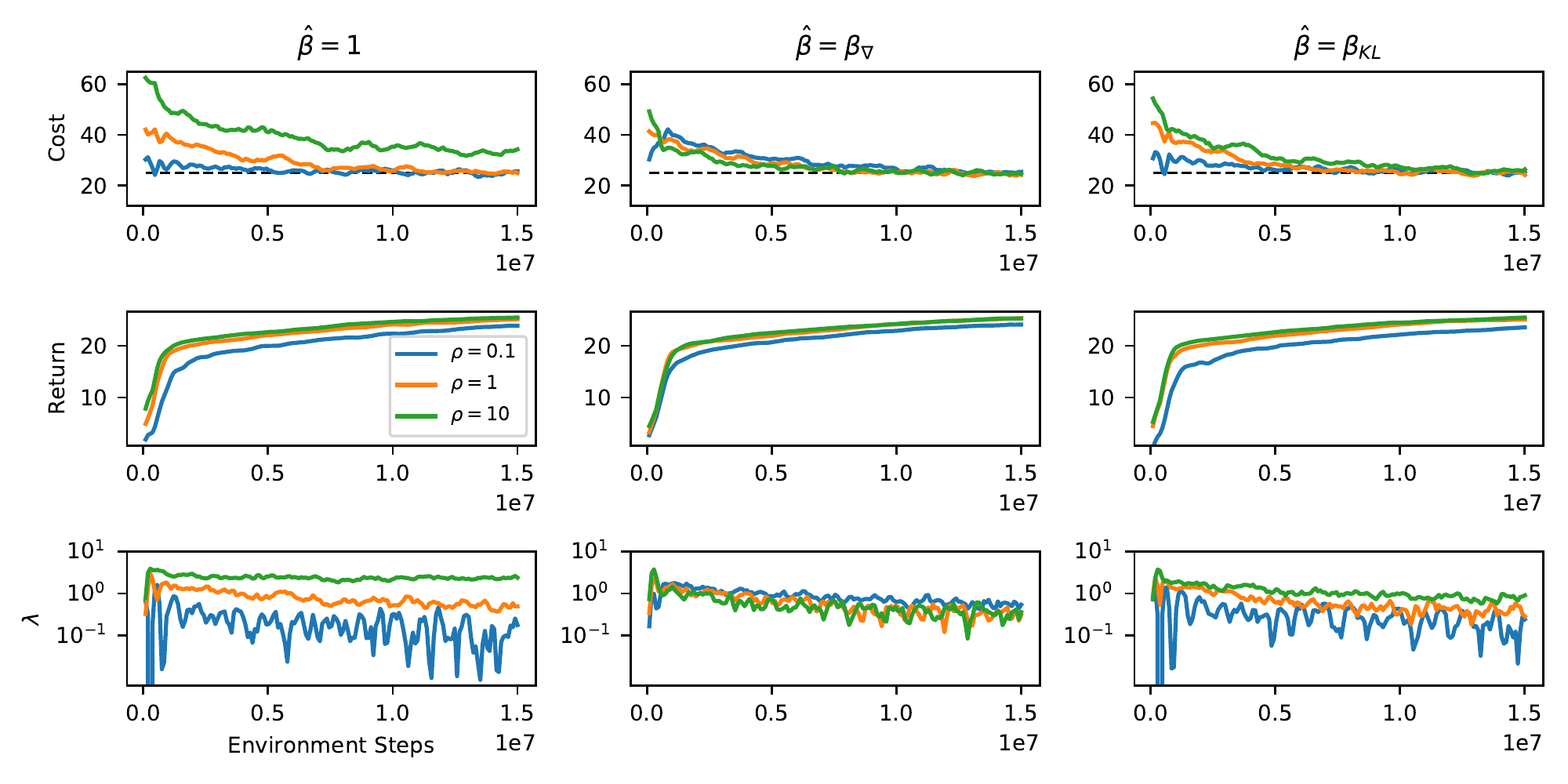}}
\underline{$K_P=10$}
\centerline{
  \includegraphics[clip,width=0.8\columnwidth]{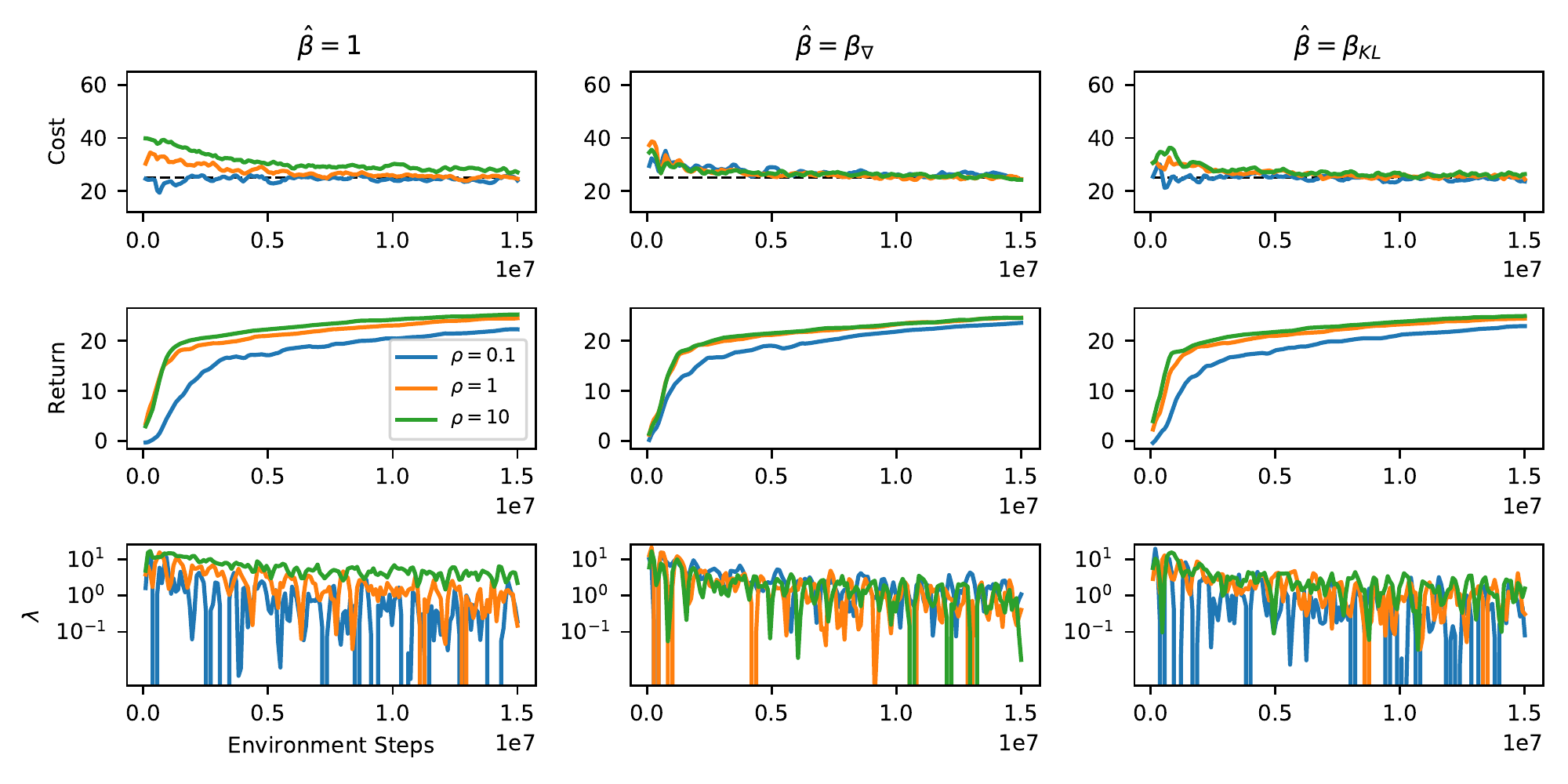}}
\caption{Reward scaling, PI-control with $K_I=0.001$, \textsc{PointGoal1}, cost-limit=25.}
\label{fig:rewardscale_pointgoal_kp}
\end{center}
\end{figure}

\begin{figure}[ht]
\begin{center}
\underline{$K_I=0.001$}
\centerline{
  \includegraphics[clip,width=0.8\columnwidth]{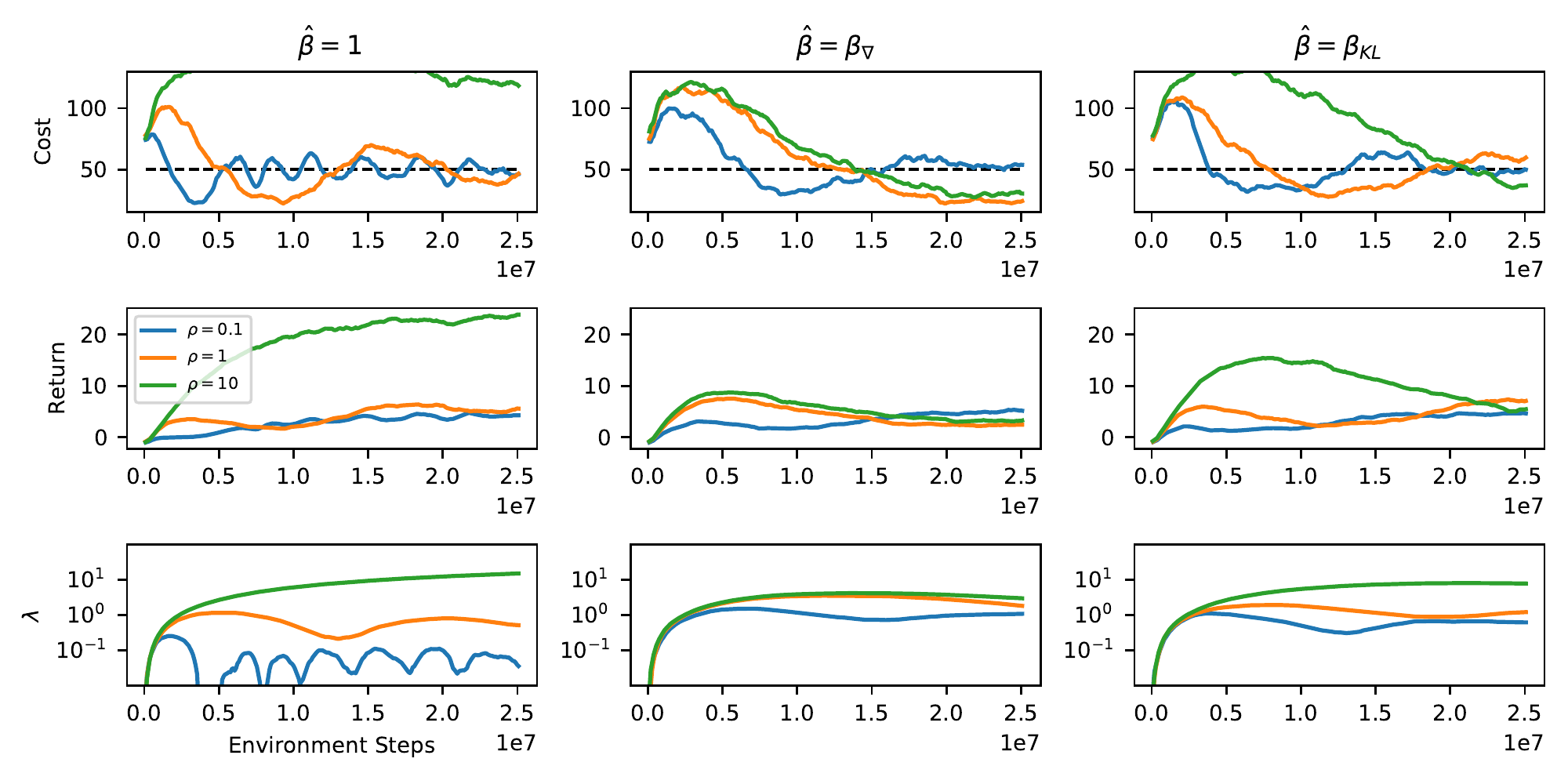}}
\underline{$K_I=0.01$}
\centerline{
  \includegraphics[clip,width=0.8\columnwidth]{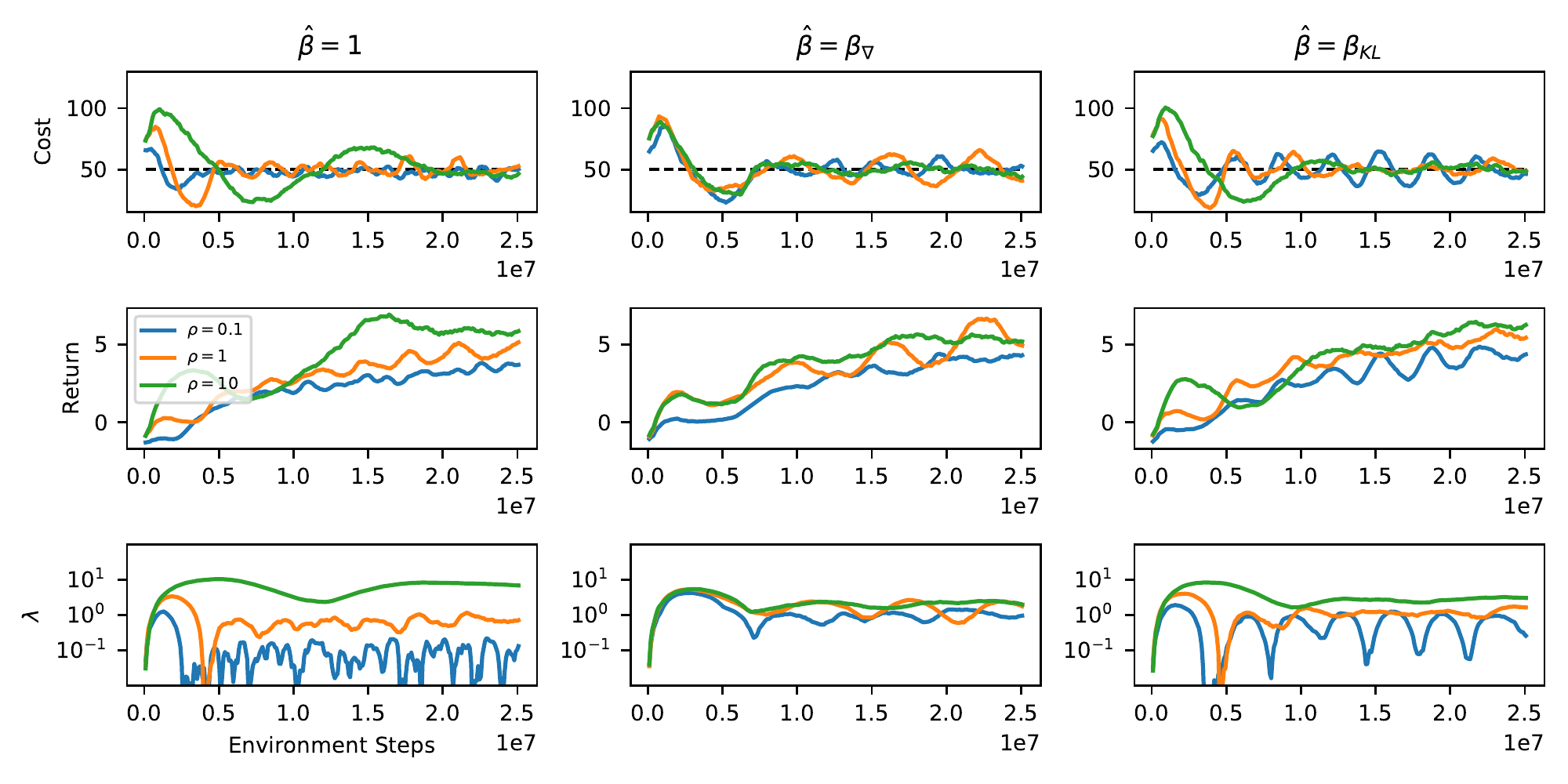}}
\underline{$K_I=0.1$}
\centerline{
  \includegraphics[clip,width=0.8\columnwidth]{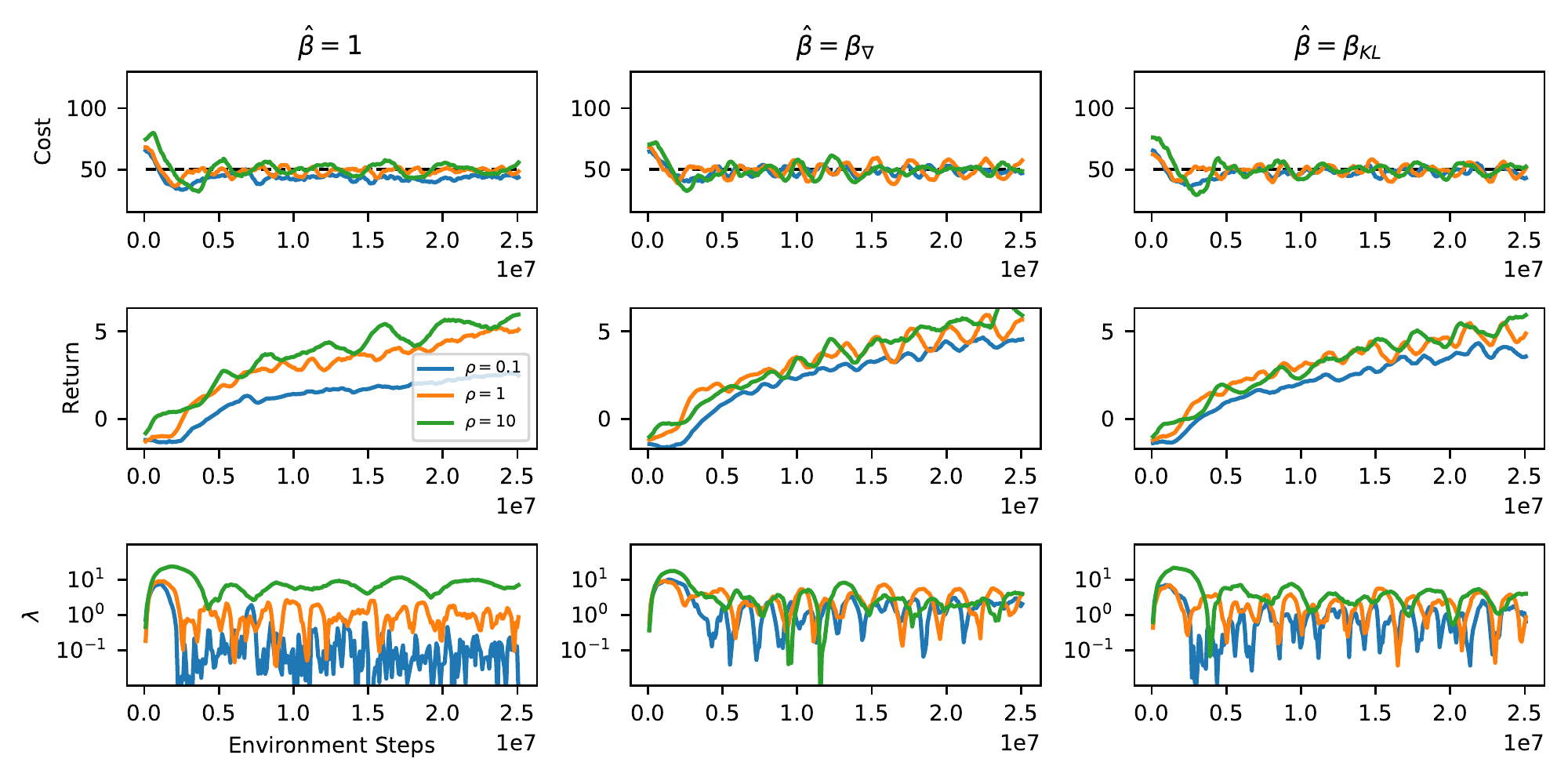}}
\caption{Reward scaling, I-control, \textsc{DoggoGoal2}, cost-limit=50.}
\label{fig:rewardscale_doggogoal}
\end{center}
\end{figure}